\documentclass{amsart}
\usepackage[centertags]{amsmath}
\usepackage{amsfonts}
\usepackage{amssymb}
\usepackage{amsthm}
\usepackage{epsfig}
\usepackage{fullpage }
\usepackage{subfigure}

\newtheorem{theorem}{Theorem}
\newtheorem{proposition}{Proposition}
\newtheorem{lemma}{Lemma}
\newtheorem{corollary}{Corollary}
\theoremstyle{definition}
\newtheorem{definition}{Definition}
\newtheorem{notation}{Notation}

\theoremstyle{remark}
\newtheorem{remark}{Remark}
\newtheorem{example}{Example}

\begin{document}

\renewcommand{\labelenumi}{(\roman{enumi})}

\newcommand{\tG}{\widetilde{G}}
\newcommand{\hG}{\widehat{G}}
\newcommand{\s}{\mathcal{S}}
\newcommand{\vx}{{\bf x}}
\newcommand{\vq}{{\bf q}}
\newcommand{\bs}{\backslash}
\newcommand{\mt}{multivariate Tutte }
\newcommand{\cH}{\s^1 (H_e)}
\newcommand{\dH}{\s^2 (H_e)}
\newcommand{\dHf}{\s^2 (H_f)}
\newcommand{\cphi}{\phi_e^{(1)}}
\newcommand{\dphi}{\phi_e^{(2)}}
\newcommand{\deta}{\eta_e^{(2)} }
\newcommand{\ceta}{\eta_e^{(1)} }
\newcommand{\doeta}{\ddot{\eta}_e^{(1)} }
\newcommand{\dteta}{\ddot{\eta}_e^{(2)} }
\newcommand{\boeta}{\bar{\eta}_e^{(1)} }
\newcommand{\calL}{{\mathcal{L}}}
\newcommand{\cG}{\check{G}}
\newcommand{\vb}{{\bf b}}
\newcommand{\BR}{Bollob\'{a}s-Riordan }
\newcommand{\al}{\alpha}
\newcommand{\be}{\beta}
\newcommand{\ga}{\gamma}
\newcommand{\de}{\delta}
\newcommand{\hs}{\hat{s}}
\newcommand{\ts}{\tilde{s}}
\newcommand{\te}{\tilde{e}}
\newcommand{\bcH}{\bar{\s}^1 (H_e)}
\newcommand{\dcH}{ \ddot{\s}^1(H_e)}
\newcommand{\bdH}{\bar{\s}^2(H_e)}
\newcommand{\ddH}{\ddot{\s}^2(H_e)}
\newcommand{\goa}{G \otimes A}
\newcommand{\goh}{G \otimes H^{p}}


\title[Expansions for polynomials of separable ribbon
graphs]{Expansions for the Bollob\'{a}s-Riordan polynomial
of separable ribbon graphs}

\author[S.~Huggett]{Stephen Huggett$^*$}
\author[I.~Moffatt]{Iain Moffatt$^\dagger$}

\begin{abstract}
We define 2-decompositions of ribbon graphs, which generalise 2-sums
and tensor products of graphs. We give formulae for the
Bollob\'{a}s-Riordan polynomial of such a 2-decomposition, and derive
the classical Brylawski formula for the Tutte polynomial of a tensor
product as a (very) special case. This study was initially motivated
from knot theory, and we include an application of our formulae to
mutation in knot diagrams.
\end{abstract}

\thanks{
${\hspace{-1ex}}^*$School of Mathematics and Statistics, University
of Plymouth, Plymouth, Devon, UK;  \\
${\hspace{.35cm}}$ \texttt{s.huggett@plymouth.ac.uk}}

\thanks{
${\hspace{-1ex}}^\dagger$
Department of Combinatorics and Optimization,  University of
Waterloo, Waterloo, Ontario, Canada; \\
${\hspace{.35cm}}$ {\em Current address:} Department of Mathematics and Statistics,  University of South Alabama, Mobile, AL 36688, USA;\\
${\hspace{.35cm}}$ \texttt{imoffatt@jaguar1.usouthal.edu}}

\date{June 28, 2007}

\maketitle

\section{Introduction}
We are interested in the decomposition of graphs and ribbon graphs
into their 2-connected components. Suppose a graph $\hG$ is
2-separable. We may regard it as arising from the 2-sums of a
collection of
graphs $\left\{{A}_e \right\}_{e \in E}$ with the graph $G=(V,E)$.
Here the
subscript $e$ plays two roles: it labels the
individual graphs in the collection $\left\{{A}_e \right\}_{e \in E}$,
and within each of these graphs it distinguishes the edge along which
the
two-sum is to be taken. The graph $G$ determines how the graphs
$\left\{{A}_e \right\}_{e \in E}$ are assembled. Strictly speaking,
the 2-sum is not well-defined on graphs without specifying which way
round the edges are to be identified: in what follows we overcome that
by referring to the vertices $u_e$ and $w_e$ at each end of the edge
$e$. 
Also, it will often be more convenient for us to work with the graphs 
$H_{e}=A_{e}\setminus\{e\}$. We will call the structure $\left(G,\{
H_{e}\}_{e\in E} \right)$ a 2-decomposition for $\hG$. There are two
important special cases. One arises when $G$ is the graph on two
vertices with two edges $e,f$ joining them. Then $\hG$ is the
conventional two-sum $A_{e}\oplus_{2}A_{f}$. The other arises when all
of the $A_{e}$ are equal to a graph $A$. Then $\hG$ is the tensor
product
$G\otimes A$ \cite{Se}. 

It is natural to seek the connection between the graph
polynomials of $\hG$ and those of $G$ and the $H_{e}$. There is a
well known result due to Brylawski~\cite{Br} which describes the
Tutte polynomial of the tensor product $T(G \otimes A)$ in terms of
those of its two factors $G$ and $A$. This result has played an
important role in the complexity theory of the Tutte polynomial
\cite{GoJe, JVW, Je}. 
Brylawski's result also plays a role in knot theory: in~\cite{Hu} the
first author used Brylawski's result to explore the relation between
the realizations of the Jones and HOMFLY polynomials as evaluations
of the Tutte polynomial of an associated graph \cite{Boll, Ja, Je,
Tr}.

Another recent example of the connection between the polynomials of
$\hG$ and those of its 2-decomposition comes from Woodall~\cite{Wo}.
In this paper he expressed the Tutte polynomial of $\hG$ in terms of
the graphs $H_{e}$ and  either the flow polynomials of subgraphs of
$G$ or the tension polynomials of contractions of $G$. This work is
related to problems on the homeomorphism classes of graphs.

\medskip

Here we are interested in generalizing the results of Brylawski in
two directions. We want to drop the condition that all of the
graphs $H_{e}$ are equal, and we also want to generalize the
formula to ribbon graphs (graphs with a cyclic ordering of each of
the incident half-edges at each vertex). This latter will
entail the study of the \BR polynomial, which is the generalization of
the Tutte polynomial to ribbon graphs.  

Brylawski's proof of the tensor product formula uses the universal
properties of the Tutte polynomial. This approach, however, only
works for
the tensor product and cannot be extended to our 2-decompositions of
graphs (although we acknowledge that Brylawski's proof does have the
advantage that it can be extended to matroids). Moreover, the
universal
properties of the \BR polynomial do not seem to be strong enough to 
support this method of proof of a Brylawski theorem for ribbon graphs
(because the basis would consist of all 1-vertex ribbon graphs). Thus
we see that a new approach is needed.

The idea behind our approach is simple. The \BR and Tutte
polynomials can be described as a sum over states, where a state is a
spanning (ribbon) subgraph. The polynomials count the number of
edges, connected components and, for the \BR polynomial, the number
of boundary cycles of the states. Let $\hG$ be as above. There
is an obvious bijection between the states of
$\hG$ and states of $\cup_{e\in E} H_{e}$. This
gives a decomposition of the states of $\hG$. We are interested in
calculating the \BR and Tutte polynomials, so we also need a way of
relating the number of connected components and boundary cycles in
the states of $H_{e}$ to those in the corresponding state of $\hG$. 
The  ribbon graph $G$ describes how each of the copies of $H_e$ are
linked together to form $\hG$, and we use the states of $G$ to
relate the states of $\hG$ and $\cup_{e\in E} H_{e}$. 

Here is a brief plan of the paper. Section 2 defines ribbon graphs
and their polynomials. In Section 3 we show how to calculate the
Tutte polynomial of $\hG$ in terms of its 2-decomposition,
generalizing the result of Brylawski. We extend our methods to ribbon
graphs and the \BR polynomial in Section 4 and show that in special
cases we can calculate the \BR polynomial of a ribbon graph from its
2-decomposition. In Section 5 we show how the \BR polynomial of any
ribbon
graph can be calculated from its 2-separation by considering
geometric ribbon graphs. We give applications of our results to the
construction of ribbon graphs with the same polynomials, and we
finish in Section 6 with an application of our work to the study of
mutations in knot diagrams.

  I.M. would like to thank Anna De Mier for helpful conversations on
the Tutte polynomial. We are very grateful for the referee's careful reading and helpful comments. 

\section{Preliminaries} \label{sec:prelim}

\subsection{Ribbon graphs and 2-decompositions} \label{ss:ribbo}

A {\em ribbon graph} $G=(V,E)$ is an orientable surface
with boundary represented as the union of $V$ closed {\em disks} and
$E$ {\em ribbons}, $I \times I$, ($I=[0,1]$ is the unit interval) such that
\begin{enumerate}
\item the discs and ribbons intersect in disjoint the line segments
$\{0,1\} \times I$;
\item each such line segment lies on the boundary of precisely one
disk and precisely one ribbon;
\item every ribbon contains exactly two such line segments.
\end{enumerate}

Ribbon graphs  arise naturally as neighbourhoods of graphs embedded in orientable surfaces. In fact it is well known that ribbon graphs are equivalent to cellularly embedded graphs in an orientable surface (for example, see \cite{GT}). Any such embedding, together with a choice of orientation of the surface, induces a cyclic ordering
of the incident half-edges at each vertex of the graph. Thus an oriented  ribbon graph  is equivalent to  a graph
(possibly with multiple edges and loops) with a fixed cyclic ordering
of the incident half-edges at each of its vertices.  We will find the latter purely combinatorial description useful on occasion.

Recall that a graph is said to be {\em n-separable} if there exists
a set of $n$ vertices whose removal disconnects the graph.
$n$-separability is a fundamental property of the structure of a
graph. By itself, this decomposition of graphs is too coarse for
ribbon
graphs as it ignores the cyclic order at the vertices and therefore
the inherent topology. For example, the graph with one
vertex and two edges $e$ and $f$ is clearly 1-separable.
However, the two choices of cyclic order of the half-edges $e,e,f,f$
and $e,f,e,f$ give rise to two distinct ribbon graphs with very
different properties, which are not captured by their 1-separable
components.

\begin{figure}
\[\epsfig{file=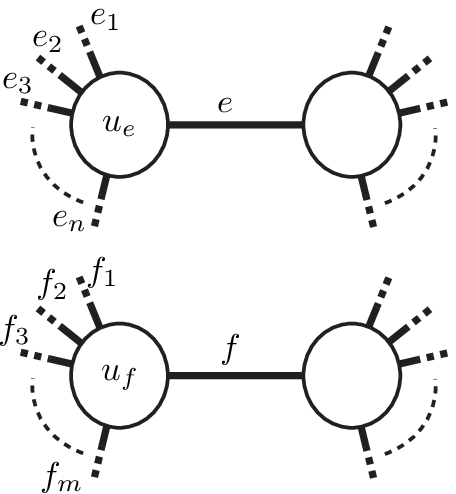, width=4cm} \quad \raisebox{2cm}{$\longrightarrow$} \quad\epsfig{file=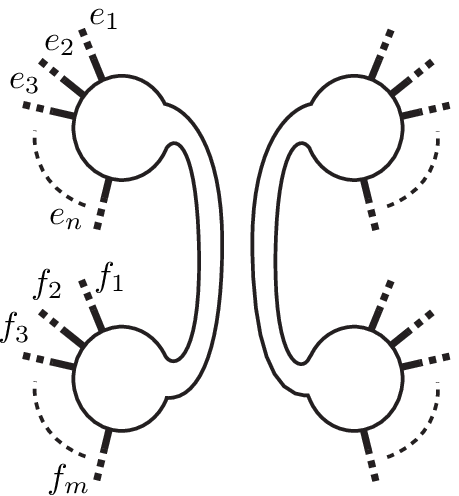, width=4cm} \quad \raisebox{2cm}{$\longrightarrow$} \quad \raisebox{1.1cm}{\epsfig{file=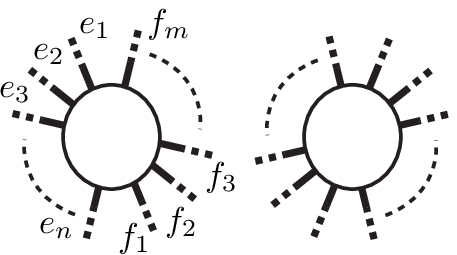, width=4cm}}\]
\caption{}
\label{fig:2sepdef1}
\end{figure}

\begin{figure}
\[\epsfig{file=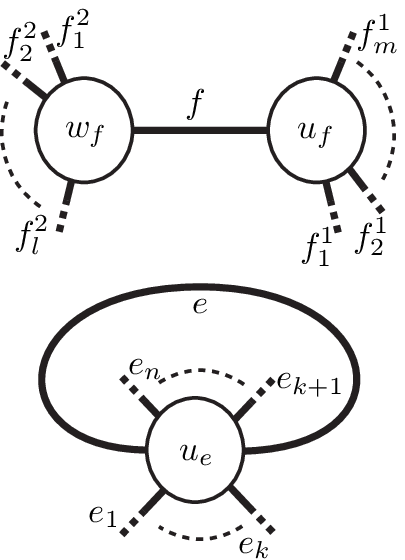, width=4cm} \quad \raisebox{2.5cm}{$\longrightarrow$} \quad\epsfig{file=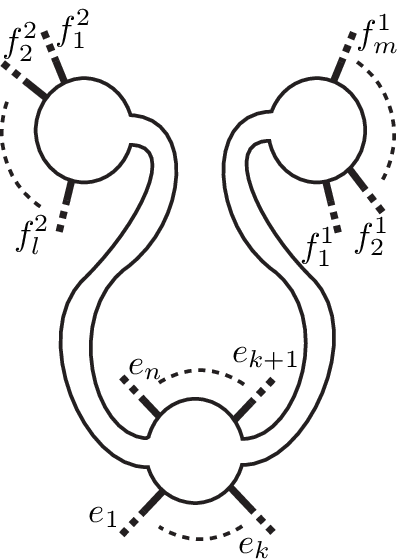, width=4cm} \quad \raisebox{2.6cm}{$\longrightarrow$} \quad \raisebox{1.2cm}{\epsfig{file=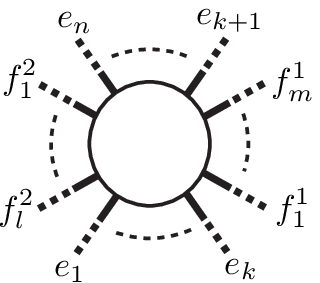, width=3.3cm}}\]
\caption{}
\label{fig:2sepdef2}
\end{figure}

However, the notion of a 2-sum is more subtle: given two graphs
$G$ and $F$ with distinguished edges $e\in E(G)$ and
$f\in E(F)$, the 2-sum $G\oplus_{2}F$ is defined by
identifying $e$ with $f$ (making a choice of which way round
to do this) and then deleting the identified edge. This works just as 
well when $G$ and $F$ are ribbon graphs: in the process of
identifying $e$ with $f$, suppose that the vertex
$u_{e}\in V(G)$ is identified with $u_{f}\in V(F)$.
We choose orientations for these ribbon graphs, and
suppose further that $e$ and $f$ are not loops and the cyclic orders
around these two vertices were 
$\{e,e_{1},\dots,e_{n}\}$ and $\{f,f_{1},\dots,f_{m}\}$. Once the
identified edge has been deleted, the cyclic order around the new
vertex will be $\{e_{1},\dots,e_{n},f_{1},\dots,f_{m}\}$. This is illustrated in figure~\ref{fig:2sepdef1}.
In the case when $e$ is a loop and $f$ is not we suppose that $u_e$
has cyclic ordering $\{e,e_{1},\dots ,e_k ,e, e_{k+1}   ,e_{n}\}$ 
and the vertices $u_f$ and $w_f$ of the edge $f=(u_f, w_f)$ have
cyclic orderings 
$\{f,f^1_{1},\dots,f^1_{m}\}$ and $\{f,f^2_{1},\dots,f^2_{l}\}$
respectively. Then once the
identified edge has been deleted, the cyclic order around the new
vertex will be $\{e_{1},\dots,e_{k},f^1_{1},\dots,f^1_{m}
,e_{k+1},\dots,e_{n},f^2_{1},\dots,f^2_{l}\}$.  This is illustrated in figure~\ref{fig:2sepdef2}.

We  would like to emphasise the fact that in the formation of the 2-sum we arbitrarily  assigned orientations to the the two summands and that the 2-sum is an orientation preserving operation. In general, the resulting ribbon graph will depend upon the choice of orientations used in the formation of the 2-sum.

\begin{definition}
Let $G=(V,E)$ be a (ribbon) graph and $\left\{A_e \right\}_{e\in E}$
be a set of (ribbon) graphs each of which has a specific
non-loop edge
distinguished. For each $e\in E$ take the 2-sum $G\oplus_{2}A_{e}$,
along the edge $e$ and the distinguished edge in $A_{e}$, to obtain
the (ribbon) graph $\hG$. For each $e\in E$ define
$H_{e}=A_{e}\setminus \{e\}$. We will call the structure $\left(G,\{
H_{e}\}_{e\in E} \right)$ a 2-decomposition for $\hG$.
The (ribbon) graph $G$ is called the {\em template}.
\end{definition}

Examples of 2-decomposition are given in Figures~\ref{fig:decomp},
\ref{fig:example}, and \ref{fig:Cydecomp}. 

Note that each of the (ribbon) graphs $H_{e}$ has two distinguished
vertices: those which were joined by the distinguished edge in
$A_{e}$.
We will denote these two distinguished vertices of each $H_e$ as $u_e$
and $w_e$ throughout this paper, and we will also use $u_e$ and $w_e$
to
denote the vertices of the edge $e$ of $G$. We will assume (without
loss
of generality) that the vertices $u_e$ and $w_e$ of $H_e$ are
distinct,
but we make no such assumption on the vertices $u_e$ and $w_e$ of
$G$, so
$e$ may be a loop. This means that the 2-sum of ribbon graphs may be
a 1-sum of  graphs.

If each of the (ribbon) graphs $H_e$ in a $2$-decomposition are equal
to a (ribbon) graph $H$, and if $u_e$ and $w_e$ lie in the same
connected component, then $\hG$ is the {\em tensor product} of $G$
with $A$, written $\hG = G\otimes A$.

Finally, if $G=C_{2}$, the 2-cycle with edges $e$ and $f$, then $\hG$
is the {\em 2-sum} $A_{e}\oplus_{2}A_{f}$.

\medskip

\subsection{Tutte and \BR polynomials} \label{ss:TandBR}Again let
$G=(V,E)$ be a (ribbon) graph. A {\em state} of $G$ is a
spanning (ribbon) subgraph $(V,E')$ where $E'\subseteq E$. We denote
the set of states by $\s(G)$. We will often abuse notation and write
``state'' when we mean ``a subset $E'\subset E$'' and vice
versa.  In particular, we will often write ``$e\in s$'' rather than ``$e\in E(s)$''.  This abuse of notation should cause no confusion.
If $s = (V,E') \in \s (G)$ then  we define $v(s):= |V|$,
$e(s):=|E'|$, $\neg e(E') := |E| - e(E')$,  $r(s):= |V|- k(s)$ and 
$n(s):=
e(s)- r(s)$, where $k(s)$  denotes the number of connected components
of $s$. Regarding a ribbon graph as a surface, we set
$\partial(s):=|\partial(s)|$, the number of its boundary components.
By a {\em planar} ribbon graph we mean one that can be regarded as a
genus zero surface.

We are interested in expansions for the Tutte polynomial of a graph
and
the \BR polynomial \cite{BR1, BR}, which is the natural
generalization of the Tutte polynomial to ribbon (or embedded) graphs.
The {\em \BR polynomial for ribbon graphs} is defined as the sum
\begin{equation}\label{eq:BRdef}
R(G;\al , \be, \ga)  = \sum_{s\in
\s(G)}(\al-1)^{r(G)-r(s)}\be^{n(s)}\ga^{k(s)-\partial(s)+n(s)}.
\end{equation}

For example, 
\[  R\left(  \raisebox{-3mm}{\epsfig{file=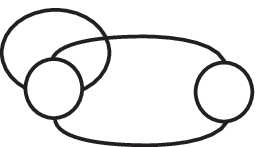, height=8mm}}   ; \alpha, \beta, \gamma \right) = \beta^2 \gamma^2 + \alpha + 2\beta +2.   \]

Observe that when  $\ga =1$ this polynomial becomes the Tutte
polynomial:
\[
R(G; x ,y-1 ,1) = T(G;x,y) :=  \sum_{s\in
\s(G)}(x-1)^{r(G)-r(s)}(y-1)^{n(s)} .
\]
This coincidence of polynomials also holds when the ribbon graph $G$
is planar. In fact the exponent of $\ga$ is twice the genus of the
state.

Before proceeding we also note that  expanding the definitions of
$r(s)$ and $n(s)$ in (\ref{eq:BRdef}) yields the following equivalent
form of the Bollob\'{a}s-Riordan, and hence Tutte, polynomials:
\begin{equation} \label{eq:BRzform}
 R(G;\al , \be, \ga)= (\al-1)^{-k(G)} (\be \ga)^{-v(G)} \sum_{s\in
\s(G)} ((\al-1)\be\ga^2)^{k(s)} (\be\ga)^{e(s)} (\ga)^{-\partial(s)}.
\end{equation}
This rewriting of the polynomials is fundamental to our approach.

\medskip

We need to consider multivariate generalizations of the Tutte and \BR
polynomials. (We will then specialize to the above polynomials where
appropriate).
The multivariate Tutte polynomial has been used extensively, and we
refer the reader to Sokal's survey article
\cite{So} for an exposition of its properties. The multivariate \BR
polynomial is the obvious extension of the multivariate Tutte
polynomial. It has been used previously \cite{CP, LM, Mo}.

The {\em multivariate \BR polynomial} of a ribbon graph $G=(V,E)$ is
\begin{equation}\label{eq:multBR}
Z(G; a, \vx ,c  ) = \sum_{s \in \s (G)} a^{k(s)} \left(\prod_{e \in
s} x_e \right) c^{\partial(s)} \; \; \in \mathbb{Z} \left[a,\{x_e
\}_{e \in E}, c \right] ,
\end{equation}
where $\vx$ denotes the set $\{x_e  \}_{e \in E}$.
The {\em multivariate Tutte polynomial} (\cite{So}) is then the
specialization
\begin{equation}\label{eq:multTutte}
Z(G; a,\vx ) := Z(G; a,\vx,1  ) = \sum_{s \in \s (G)} a^{k(s)}
\left(\prod_{e \in s} x_e \right).
\end{equation}

When we set all of the variables $x_e$ equal to $b$, say, then we
will denote the polynomial $Z(G; a,\vx,c )|_{x_e=b}$ simply by $Z(G;
a,b,c )$. Notice that by equation~(\ref{eq:BRzform}),
\begin{equation}\label{eq:RZ}
R(G;\al , \be, \ga)= (\al-1)^{-k(G)} (\be \ga)^{-v(G)}
Z\left(G;(\al-1)\be\ga^2, \, \be\ga, \, \ga^{-1}\right),
\end{equation}
so that the specialization is equivalent to the  \BR polynomial. A
similar relation holds for (\ref{eq:multTutte}) and the Tutte
polynomial.

When the choice of variables is clear from the context we will just
write $Z(G)$ instead of $Z(G; a, \vx ,c)$ or its specializations
(such as $x_e = b$ or $c=1$).


\begin{notation}\label{notation1}
Henceforth $\hG$ will always denote a (ribbon) graph with
the 2-decomposition $\left(G,\{H_e\}_{e\in E}\right)$. Here $H_e$
denotes the graph $A_e$ with its distinguished edge (joining the
vertices $u_e$ and $w_e$) deleted. We will often only need to pick out
this distinguished edge in $Z(A_{e};a,\vx,c)$, so we denote by $\vb$
the specialization of $\vx$ which leaves $x_{e}$ as it is but for all
$d\,(\neq e)\in E(A_{e})$ sets $x_d=b$.
In a ribbon graph the distinguished edge will be between the marked
points $m_e$ and $n_e$ described later in Subsection~\ref{ss:CyDe}
and Figure~\ref{fig:dmark}.
\end{notation}

\subsection{Deletion and contraction for ribbon graphs}
Later we will find it convenient to express our results in terms of
the deletion and contraction of a ribbon in a ribbon graph. Here, we
use
the surface description of a ribbon graph. 
The deletion of a ribbon always makes sense, but we need to be
more careful when we contract a ribbon.

Let $e$ be a ribbon of $F$.
First suppose that $e$ is not a loop. Then $e$ is a ribbon between
two distinct disks $u$ and $w$ of $F$. Suppose that the cyclic order
of incident half-edges is $e_{u1}, \ldots , e_{un}, e$ at $u$, and
$e, e_{w1}, \ldots , e_{wn}$ at $w$. Then $F/e$ is the ribbon graph
obtained by replacing the vertices $u$ and $v$ with a single vertex
with cyclically ordered incident edges $e_{u1}, \ldots , e_{un},
e_{w1}, \ldots , e_{wn}$.

Next suppose that  $e$ is a loop. To describe the contraction of a
loop we need a generalization of ribbon graphs, by allowing
identifications of the boundary of any disk in a ribbon graph. We
define a {\em ribbon surface} to be a surface with boundary
represented
as the union of $|V|$ closed surfaces with boundary, called {\em
nodes}, and $|E|$ {\em ribbons}, $I \times I$, such that
\begin{enumerate}
\item the nodes and ribbons intersect in disjoint line segments
$\{0,1\} \times I$;
\item each such line segment lies on the boundary of precisely one
node and precisely one ribbon;
\item  every ribbon contains exactly two such line segments.
\end{enumerate}
The contraction of  a loop $e$ is then defined by identifying the 
 endpoints $\{0\}\times I$ and $\{1\} \times I$
of $e$ and then deleting $e$.

Note that if all of the nodes are disks, then an orientable ribbon
surface is exactly a ribbon graph. We carry over all of our notation
for ribbon graphs, except that by $v(F)$ we mean the number of nodes. 
With this inherited notation the \BR polynomial of a ribbon surface
makes sense. Moreover the following deletion-contraction  relation
now holds for any edge $e$:
\begin{equation}\label{eq:delcont}
Z(F; a, \vx ,c  ) = Z(F\backslash e ; a, \vx ,c  ) + x_e Z(F/ e ; a,
\vx ,c  ) .
\end{equation}
Notice that this generalizes the deletion-contraction relations for
the
\BR polynomial given in \cite{BR} for ordinary edges and trivial
loops.

We shall abuse notation and extend to ribbon surfaces without comment
when the need arises. It is possible to avoid the generalization to
ribbon surfaces, but it is much more convenient and concise to be
able to use a deletion-contraction relation for loops.

\subsection{Geometric ribbon graphs}\label{sec:geo}
A  geometric ribbon graph is defined just as in the surface
realization of a ribbon graph but without the condition of
orientability, so we allow the ribbons between the disks to
be twisted any number of times.

We carry over all of the notation for ribbon graphs, in particular
that of Subsection~\ref{ss:ribbo}, to geometric ribbon graphs. The
one additional piece of information we need from a state of $F$ is
whether it is
orientable or not. Let $s \in \s (F)$; we set $t(s)=0$ if $s$ is
orientable and  $t(s)=1$ if $s$ is non-orientable.

The {\em \BR polynomial for geometric ribbon graphs} (see \cite{BR})
is defined as follows.
\begin{equation*}\label{eq:BRgeodef}
R(F;\al , \be, \ga, \de)  = \sum_{s\in
\s(F)}(\al-1)^{r(F)-r(s)}\be^{n(s)}\ga^{k(s)-\partial(s)+n(s)}
\de^{t(s)}  \; \; \in \mathbb{Z} \left[\al,\be,\ga ,\de
\right]/(\de^2-\de) .
\end{equation*}
Observe that when  $\de =1$ or $F$ is orientable this polynomial
coincides with the \BR polynomial~(\ref{eq:BRdef}). 

As with the ribbon graph polynomial, we need to use a multivariate
version of this polynomial.
The {\em multivariate \BR polynomial} of a geometric ribbon graph
$F=(V,E)$ is
\begin{equation*}\label{eq:multBRgeo}
Z(F; a, \vx ,c,d) = \sum_{s \in \s (F)} a^{k(s)} \left(\prod_{e \in
s} x_e \right) c^{\partial(s)} d^{t(s)} \; \; \in \mathbb{Z}
\left[a,\{x_e
\}_{e \in E},c ,d \right]/(d^2-d) ,
\end{equation*}
where $\vx$ denotes the set $\{x_e  \}_{e \in E}$.

When we set all of the variables $x_e$ equal to $b$, say, then we
will denote the polynomial  by $Z(G;a,b,c,d )$. We then have
\[
R(G;\al , \be, \ga,\de)= (\al-1)^{-k(G)} (\be \ga)^{-v(G)}
Z\left(G;(\al-1)\be\ga^2, \, \be\ga, \, \ga^{-1}, \, \de \right).
\]

\begin{remark}
Although we have favoured a surface description of geometric ribbon
graphs here, everything could have been defined in terms of
 graphs with a fixed cyclic ordering of the incident half-edges
at each vertex by adding a $+$ or $-$ sign on each edge to record the parity of the number of half-twists. We refer the reader to Bollob\'{a}s and Riordan's
article \cite{BR} for details. 
\end{remark}


\section{The Tutte polynomial}\label{sec:tutte}

\subsection{The Decomposition} \label{ss:TuDe}                    
Suppose we are given a state $\hs \in \s (\hG)$ and a 2-decomposition
$\left(  G , \{ H_e \}_{e\in E}   \right)$
of  $\hG$. The state $\hs$ is uniquely determined by a set of edges,
but each of these edges also belongs to a graph in the set $\{ H_e
\}_{e\in E}$. Therefore the state $\hs$ uniquely determines a set of
states
$\{ s_e \in \s (H_e)\}_{e\in E}$.

Now by the definition of 2-decomposition, each graph $H_e$ has two
distinguished vertices $u_e$ and $w_e$ which are also vertices of
$G$. For each $e$ we partition the states of $\s (H_e)$ into two
subsets: $\cH$ consists of all states in $\s (H_e)$ in
which $u_e$ and $w_e$ lie in the same connected component, and
$\dH$ consists of all states in $\s (H_e)$ in which $u_e$
and $w_e$ lie in different connected components. 

\begin{figure}
\[\epsfig{file=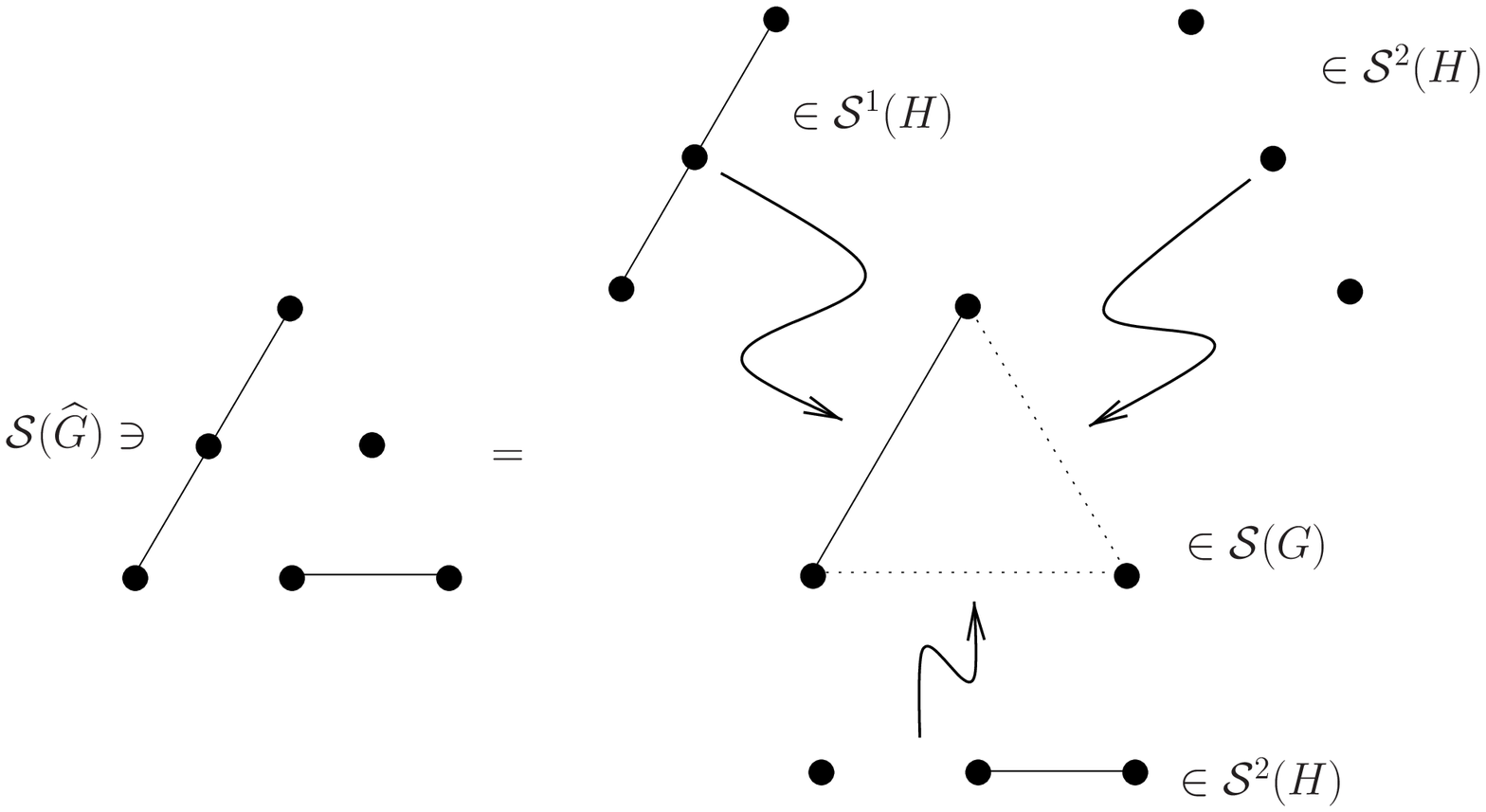, width=12cm}\]
\caption{}
\label{fig:decomp}
\end{figure}

We will use the partition of each set $\s (H_e)$ to construct a state
of $G$ from the states $\{ s_e\}$ determined by $\hs \in \s (\hG)$ as
above. To do this, start with the graph $G$ and then remove an edge
$e=(u_e, w_e)$ if and only if the state $s_e$ determined by $\hs$
lies in the set $\dH$.
An example is shown in Figure~\ref{fig:decomp} where $\hG=C_3\otimes
C_3$.

If we replace the edges $e$ in a state $s \in \s(G)$ with elements of
$\cH$, and the edges $f$ which are not in $s$ by elements of
$\dHf$,
we obtain a state $\hs \in \s (\hG)$. Moreover, each state of $\hG$ is
uniquely obtained in this way: we could start with $\hs \in \s (\hG)$
and from it determine an element of $\cH \cup \dH$ for each $e \in E$,
and an element of $\s (G)$.

\begin{lemma}\label{lem:con}
If a state $\hs \in \hG$ is decomposed into states
$s \in \s(G)$, $s_e   \in \cH \cup \dH$, $e \in E$, then
\[ k(\hs) =   \sum_{e \in E} k(s_e) -   |\{ s_e \in \cH \} | - 2| \{
s_e \in \dH \} | + k(s)  . \]
\end{lemma}

\begin{proof}
Each component of $s$ corresponds to a component of $\hs$, but $\hs$
has extra
components, arising from states $s_e \in \cH \cup \dH$ which have
components not containing any vertices of $G$. For each $e \in E$ for 
which $s_e \in \cH$ there are $k(s_{e})-1$ of these extra components, 
while if $s_e \in \dH$ there are $k(s_{e})-2$ of them.
    
Therefore
\[ k(\hs) =   \sum_{e \in E} k(s_e) -   |\{ s_e \in \cH \} | - 2| \{
s_e \in \dH \} | + k(s)  \]
as required.
\end{proof}

\subsection{An expansion for the Tutte polynomial}\label{ss:TuExp}
We work with the Tutte polynomial in the form
\begin{equation}
Z(G;a,b)=\left(\frac{a}{b}\right)^{k(G)}b^{v(G)}T(G;\frac{a}{b}+1,b+1),
\end{equation}
where $Z(G;a,b)=  \sum_{s\in \mathcal{S}(G)} a^{k(s)}b^{e(s)} $. 
\begin{lemma}\label{lem:tutteOGS}
Let $\left(G , \{ H_e \}_{e\in E}\right)$ be a 2-decomposition of 
$\hG$. Then
\begin{equation*}
Z(\hG ;a,b) =\sum_{s \in \s(G)} a^{k(s)} \left( \prod_{e \in s} \cphi
\right)  \left(
\prod_{e \notin s} \dphi \right)
\end{equation*}
where 
\[\begin{array}{l}
\cphi (a,b) :=  \sum_{s \in \cH} a^{k(s)-1}
b^{e(s)}, \\
\dphi (a,b)  := \sum_{s \in \dH} a^{k(s)-2}
b^{e(s)}.
\end{array}\]
\end{lemma}

\begin{proof}
    
If a state $\hs$ of $\hG$ is decomposed into states $s\in\s(G)$,
$s_e\in\cH$ and $t_f\in \mathcal{S}^2(H_f)$, then by Lemma~\ref{lem:con} we have 
\begin{equation*}
k(\hs) = k(s) + \sum (k(s_e)-1) + \sum (k(t_f)-2); 
\end{equation*}
and clearly
\begin{equation*}
e(\hs) =  \sum e(s_e) + \sum e(t_f). 
\end{equation*}

Now pick a state $s\in\s (G)$. Then for each edge $e\in s$ pick a
state
$s_{e}\in \s^{1}(H_{e})$ and for each edge $f\notin s$ pick a state
$t_{f}\in \s^{2}(H_{f})$. Each state in $\s(\hG)$ is uniquely
obtained in this way, and we obtain the term
\begin{eqnarray*}
    & &a^{k(s)}\Pi_{e\in s}a^{k(s_{e})-1}b^{e(s_{e})}\Pi_{f\notin
    s}a^{k(t_{f})-2}b^{e(t_{f})}\\
    &=&a^{k(s)+\sum (k(s_{e})-1)+\sum (k(t_{f})-2)}b^{\sum
    e(s_{e})+\sum e(t_{f})}\\
    &=&a^{k(\hs)}b^{e(\hs)}
\end{eqnarray*}
Summing over all of the states of $\hG$ then gives the result.

\end{proof}

\begin{remark}
This proof is really just the composition and sum lemmas for certain
ordinary generating series. Indeed much of this work may be expressed 
in terms of ordinary generating series.
\end{remark}

\begin{example}
As a simple example of Lemma~\ref{lem:tutteOGS}, consider the
2-decomposition $\left(G,H_{f},H_{g}\right)$ of the graph $\hG$ shown
in Figure~\ref{fig:example}.
\begin{figure}
\[ \epsfig{file=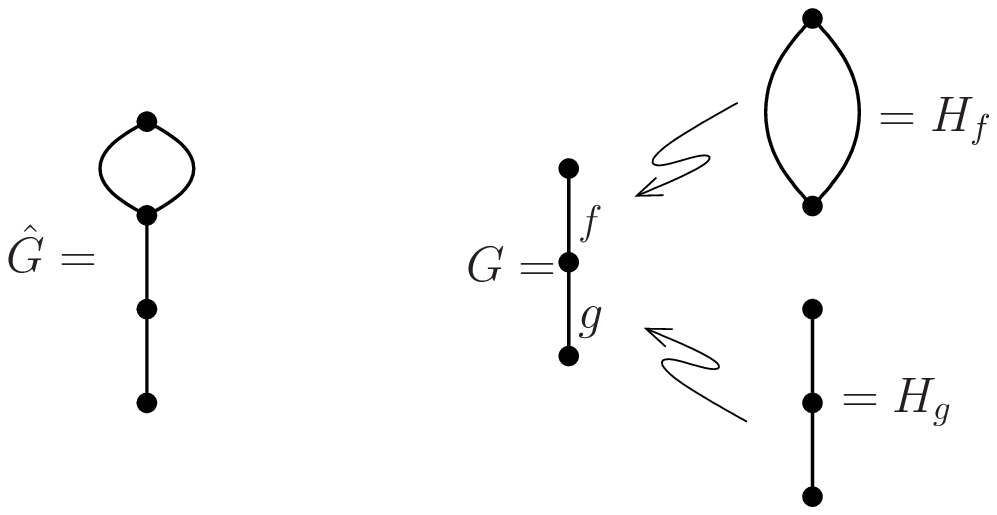, width=7cm}   \]
\caption{}
\label{fig:example}
\end{figure}
We have
\begin{gather*}
\sum_{s\in\s(G)} a^{k(s)} \left(\prod_{e\in s} \phi_{e}^{(1)} \right) \left(
\prod_{e\notin s} \phi_{e}^{(2)} \right)
=a\phi_{f}^{(1)}\phi_{g}^{(1)} + a^2\phi_{f}^{(1)}\phi_{g}^{(2)}  +a^2\phi_{f}^{(2)}\phi_{g}^{(1)}  + a^3\phi_{f}^{(2)}\phi_{g}^{(2)},\\
\phi_{f}^{(1)}(a, b) = b^2 + 2b,  
 \quad
\phi_{f}^{(2)}(a, b) = 1,  
 \quad
\phi_{g}^{(1)}(a, b) = b^2,   
 \quad
\phi_{g}^{(2)}(a, b) = 2b + a.
\end{gather*}
Therefore
\begin{eqnarray*}
\sum_{s \in \s(G)} a^{k(s)} \left( \prod_{e \in s} \cphi
\right)  \left(
\prod_{e \notin s} \dphi \right)
&=&  a(b^2+2b)b^2+a^2(b^2+2b)(2b+a)+a^2(1)b^2+a^3(1)(2b+a) \\
&=&
ab^4 + 2ab^3 + a^2b^2 + 2a^2b^3 + a^3b^2 + 4a^2b^2 + 2a^3b + 2a^3b +
a^4 \\
&=& \sum_{\hs\in\s(\hG)}a^{k(\hs)}b^{e(\hs)}
\end{eqnarray*}
as required.
\end{example}

Recalling the definition of the multivariate Tutte polynomial
(\ref{eq:multTutte}), we see that Lemma~\ref{lem:tutteOGS} can be
written as
\begin{equation}\label{eq:TuOGS}
Z (\hG;a,b) =
\left(\prod_{e\in E}\dphi\right)
Z\left(G; a, \left\{  \cphi / \dphi \right\}_{e\in E}   \right).
\end{equation}

It remains to express $\cphi$ and $\dphi$ in terms of the Tutte
polynomials of certain graphs.

\medskip

Consider the graph $A_e$, defined in Notation~\ref{notation1}.
The states of $A_e$ occur in pairs $s$ and $s\cup e$ where $e\notin
s$. Suppose we have a state $s$ that does not contain $e$. Then $s$
contributes a term $a^{k(s)}b^{e(s)}$ to the polynomial
$Z(A_e;a,\vb)$, where $\vb$ is as in Notation~\ref{notation1}. Now
$s\cup e$ has one additional edge with variable $x_e$ and it will 
either have $k$ connected components, if $u_e$ and $w_e$ are contained
in the same connected component of the state $s$, or it will have
$k-1$
connected components, if $u_e$ and $w_e$ lie in different components
of
$s$. Thus we may write
\begin{align*}
Z(A_e;a,\vb) &= (1+x_{e})\sum_{s \in \cH} a^{k(s)} b^{e(s)}+
(1+a^{-1}x_{e}) \sum_{s \in \dH} a^{k(s)} b^{e(s)} \\
&=(1+x_{e})a \cphi + (1+a^{-1}x_{e})a^{2} \dphi.
\end{align*}
By separating the terms containing $x_e$ we can also write
\[
Z(A_e;a,\vb) = Z(H_e;a,b) + x_{e}Z(A_e/e;a,b).
\]
So we can determine $\cphi$ and $\dphi$ as the unique solution to the
equations
\begin{align*}
a\cphi + a^{2}\dphi &=  Z(H_e;a,b) \\
a^{}\cphi + a^{} \dphi &= Z( A_e/e;a,b).
\end{align*}

Collecting this together we have:
\begin{theorem}
Let $\hG$ be a graph which has been obtained from the graph $G=(V,E)$
by 
taking successive two-sums along each edge $e=(u_e, w_e)\in E$ with
graphs $A_e$. In addition, let $H_e$ be the graph $A_e$ with the edge
$e=(u_e, w_e)$ deleted. Then
\[
Z(\hG;a,b) = \left( \prod_{e \in E}g_e \right)
Z\left(G;a,\left\{  f_e / g_e \right\}_{e\in E} \right).
\]
where $f_e$ and $g_e$ are the solutions to
\begin{align*}
a(f_e + ag_e) &=  Z(H_e;a,b) \\
a(f_e + g_e) &= Z(A_e /e;a,b).
\end{align*}
\hfill$\Box$
\end{theorem}

\begin{remark}
We note that by making a few trivial changes to the argument, the
2-variable polynomial $Z(\widehat{G};a,b)$ can be replaced by the multivariate 
Tutte polynomial $Z(\widehat{G};a,\vx)$ in the theorem. 
\end{remark}

Since we were considering 2-decompositions of graphs where each $H_e$
was allowed to be distinct we were forced into considering the
multivariate Tutte polynomial: we needed some way of recording which
$H_e$ went where. However, if we insist that all of the graphs $H_e$
are equal to $H$ say, then we do not need the multivariate Tutte
polynomial and we have
\begin{corollary}\label{cor:Ztensor}
\[
Z(G \otimes A;a,b) = (g)^{e(G)}Z\left( a,  f / g \right).
\]
where $f$ and $g$ are the solutions to
\begin{align*}
a(f + ag) &= Z(H;a,b) \\
a(f + g) &= Z(A/e;a,b).
\end{align*}
\end{corollary}

\begin{example}\label{ex:1}
Let $A = C_3$, the 3-cycle. Then $g=a+2b$ and $f= b^2$ and we have
\[Z(G \otimes C_3;a,b)=(a+2b)^{e(G)} Z(G;a,b^2/(a+2b)).\]
\end{example}

Since $Z(F)$ is just a rewriting of the Tutte polynomial we can use
the above corollary to recover Brylawski's theorem:
\begin{corollary}\label{cor:bry}
\[
T(G \otimes A;x,y) = (h)^{n(G)} (h^{\prime})^{r(G)}
T\left(G;T(H ;x,y)/h^{\prime},T(A/e;x,y)/h \right),
\]
where $h$ and $h^{\prime}$ are the unique solutions to
\begin{align*}
(x-1)h + h^{\prime} &= T(H;x,y) \\
h + (y-1)h^{\prime} &= T(A/e;x,y).
\end{align*}
\end{corollary}

\begin{proof}
First of all note the identity
\begin{equation} \label{eq:TutteZ}
    T(F;x,y)=(x-1)^{-k(F)}(y-1)^{-v(F)}Z(F;(x-1)(y-1),y-1),
\end{equation}
for any graph $F$. This rewriting allows us to apply
Corollary~\ref{cor:Ztensor} to the Tutte polynomial $T(G \otimes
A;x,y)$ as follows.
\begin{equation*} \begin{split}
T(G \otimes A;x,y) & = 
(x-1)^{-k(\goa)} (y-1)^{-v(\goa)} Z(\goa; (x-1)(y-1), y-1)
\\  &
=  (x-1)^{-k(\goa)} (y-1)^{-v(\goa)} g^{e(G)} Z(G; (x-1)(y-1), f/g)
\\  &
=   (x-1)^{-k(\goa)} (y-1)^{-v(\goa)} g^{e(G)} (f/g)^{v(G)} (
(x-1)(y-1)g/f)^{k(G)} 
\\  & \hspace{8cm}  T\left(G;  \frac{(x-1)(y-1)g+f}{f},\frac{f+g}{g}
\right) ,
\end{split}\end{equation*}
where $f$ and $g$ are the solutions to the equations of
Corollary~\ref{cor:Ztensor}.

Now we set
\[
h:=(x-1)^{-k(H)+1}(y-1)^{-v(H)+2} g, \quad h^{\prime}:=
(x-1)^{-k(H)+1}(y-1)^{-v(H)+1} f.
\]
Our definition of the tensor product has the effect of making
$k(H)=k(A/e)$. Also, $v(A/e)=v(H)-1$, and so by (\ref{eq:TutteZ})
the linear equations of Corollary~\ref{cor:Ztensor} become
\begin{align*}
(x-1)h + h^{\prime} &= T(H;x,y) \\
h + (y-1)h^{\prime} &= T(A/e;x,y)
\end{align*}
as required. The expression for the Tutte polynomial above becomes
\begin{multline*}
T(G \otimes A;x,y) = 
 (x-1)^{-k(\goa)} (y-1)^{-v(\goa)}  
 (x-1)^{(k(H)-1)e(G)} (y-1)^{(v(H)-2)e(G)}  \\ h^{e(G)}  
 \left(\frac{(y-1)h^{\prime}}{h}\right)^{v(G)} \left(
\frac{(x-1)h}{h^\prime}\right)^{k(G)} T\left(G;
\frac{(x-1)(y-1)g+f}{f},\frac{f+g}{g} \right).
\end{multline*}
It is easily seen that $v(\goa)=(v(H)-2)e(G)+v(G)$ and
$k(\goa)=(k(H)-1)e(G)+k(G) $, and the result follows by cancelling
terms.
\end{proof}

\begin{remark}
In the definition of the tensor product of two graphs we insisted
that the vertices $u$ and $w$ of the distinguished edge $e$ of $A$ lay
in the same connected component of $H$. This restriction was not used
in the proof of Corollary~\ref{cor:Ztensor}, but it was needed in
Corollary~\ref{cor:bry}. However if we remove this condition, then we
can improve Corollary~\ref{cor:Ztensor} as follows:\[
T(G \otimes A;x,y) = (h)^{n(G)} (h^{\prime})^{r(G)}
T\left(G;T(H;x,y)/h^{\prime},T(A/e;x,y)/((x-1)^{\varepsilon}h)\right),
\]
where $h$ and $h^{\prime}$ are the unique solutions to
\begin{align*}
(x-1)h + h^{\prime} &= T(H;x,y) \\
(x-1)^{\varepsilon} \left(h + (y-1)h^{\prime}\right) &= T(A/e;x,y),
\end{align*}
and $\varepsilon := k(H)-k(A/e)$. 
Similar statements will hold for Corollary~\ref{cor:planbryII}.
\end{remark}

\section{The \BR polynomial I:  embedding in a
neighbourhood}\label{sec:plan}

We will begin our study of the \BR polynomial of 2-decomposable
ribbon graphs with a particularly pleasing special case.

Throughout this section $\hG$ will denote an embedded graph with a
2-decomposition $\left( G , \{ H_e \}_{e\in E} \right)$, where $G$ is
embedded. We will further assume  that each graph $H_e$ is embedded
in a neighbourhood
of the edge $e$ of the embedded graph $G$ in the formation of $\hG$. 
This restriction on the embedding of $\hG$ imposes a very strong
connection between the topology of $\hG$ and $G$ which we will use to
our advantage.

\subsection{Decomposing the states} \label{ss:CyDe}

We proceed as in the case of the Tutte polynomial. 
A state $\hs \in \s (\hG)$  uniquely determines and  is uniquely
determined by
a set of
states $\{ s_e \in \s (H_e)\}_{e\in E}$.
Again each graph $H_e$ has two
distinguished vertices $u_e$ and $w_e$ which are also vertices of
$G$. For each $e \in E $ we partition the states of $\s (H_e)$ into
two
subsets: $\cH$, which consists of all states in $\s (H_e)$ in
which $u_e$ and $w_e$ lie in the same connected component; and
$\dH$, which consists of all states in $\s (H_e)$ in which $u_e$
and $w_e$ lie in different connected components. 

Given a state $\hs \in \s (\hG)$, which determines a unique set of
states $\{ s_e \in \s (H_e)\}_{e\in E}$, construct a state of $G$ by
removing an edge $e=(u_e, w_e)$ from $G$
 if and only if the state $s_e$ belongs to the partition  $\dH$. 
This is exactly the construction used in Section~\ref{sec:tutte} to
study the Tutte polynomial. 

Conversely, given a state $s \in \s(G)$ and a copy of the template
$G$, if we replace the edges of $G$ which are also in the state $s$
with elements of $\cH$, and the other edges of $G$ with elements of
$\dH$, we obtain a state in $\s (\hG)$. Clearly, each state of $\hG$
is uniquely obtained in this way.

\begin{lemma}\label{lem:brcon}
If a state $\hs$ of the embedded graph $\hG$ is decomposed into states
$s \in \s(G)$, $s_e   \in \cH \cup \dH$, $e \in E$ by the
decomposition above, then
\[ k(\hs) =   \sum_{e \in E} k(s_e) -   |\{ s_e \in \cH \} | - 2| \{
s_e \in \dH \} | + k(s)   \]
and
\[ \partial (\hs) =   \sum_{e \in E} \partial(s_e) -   |\{ s_e \in
\cH \} | - 2| \{
s_e \in \dH \} | + \partial(s) ,  \]
where $\partial$ counts the boundary components of the associated
ribbon graphs.
\end{lemma}

\begin{figure}
\[\epsfig{file=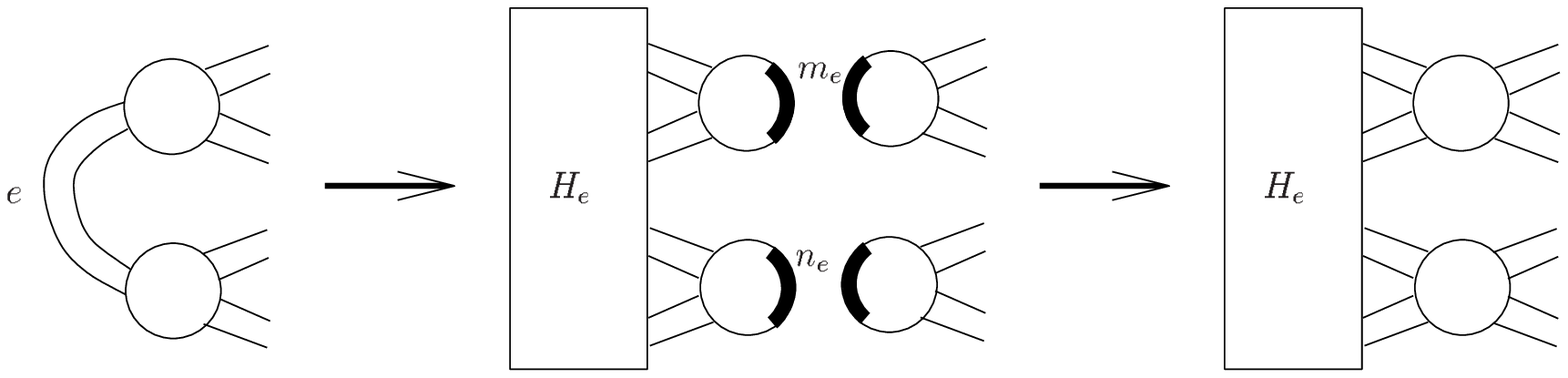, width=10cm}\]
\caption{}
\label{fig:dmark}
\end{figure}

To prove this lemma we introduce some notation, which is also useful
later.
The ribbon graph $\hG$ is obtained from $G$ by replacing each ribbon
$e$ with $H_e$. 
The two ends $\{0,1\}\times I$ of a ribbon $e$ induce two arcs on the
incident discs $u_e$ and $w_e$  ($u_e$ and
$w_e$ may be the same vertex). Denote these two arcs by $m_e$ and
$n_e$. 
We may then view the replacement of the ribbon $e$ with $H_e$ as the
operation which identifies an arc on the disc $u_e$ of $H_e$ with
one of the arcs $n_e$ or $m_e$, and an arc on $w_e$ of $H_e$ with the
other. We  also
denote these two arcs in $H_e$ by  $n_e$ and $m_e$ according to their
identification. (Note that by definition the discs $u_e$ and
$w_e$ in $H_e$ are distinct.)  
An example is shown in Figure~\ref{fig:dmark}.

When the vertices $u_e$ and $w_e$ of $H_e$ and $G$ are identified in
the formation of $\hG$, the boundary components and the connected
components containing the marked points $m_e$ are merged, and the
boundary and connected  components  containing the marked points
$n_e$ are merged. 

Notice that in the ribbon graph $H_e$, either $m_e$ and $n_e$ belong
to the same boundary component, in which case they also belong to the
same connected component; or   
$m_e$ and $n_e$ belong to distinct boundary components and may or
may not belong to the same connected component. For example, in the
ribbon graph $T-e$ of Figure~\ref{fig:gtilde}, $m_e$ and $n_e$ will
belong to different boundary components but the same connected
component. In this section we insist that each ribbon graph $H_e$
embeds into a neighbourhood of the edge $e$. This means that $H_e$ and
$A_e$ are planar. We then have 
\[  2k(A_e)-\partial (A_e)+e(A_e)-v(A_e) = 0
=2k(H_e)-\partial (H_e)+e(H_e)-v(H_e).   \] 
This gives the relation
\[2k(A_e)-\partial (A_e) = 2k(H_e)-\partial (H_e)-1.\]
Then if $m_e$ and $n_e$ belong to distinct boundary components in
$H_e$, we have  
$\partial (A_e)=\partial (H_e)-1$. 
Therefore $k(A_e)=k(H_e)-1$ so $m_e$ and $n_e$ belong to distinct
connected components.

This tells us that (under the embedding condition used in this
section) the marked points $m_e$ and $n_e$ in $H_e$ belong to the
same boundary component if and only if they belong to the same
connected component.

\begin{proof}[Proof of Lemma~\ref{lem:brcon}.]
The connectivity relation follows from Lemma~\ref{lem:con}.

As for the second identity, each boundary component of $s \in \s (G)$
corresponds to a boundary component of $\hs \in \s (\hG)$.  To see why this is, with reference to  
Figure~\ref{fig:dmark}, let $a_e$ and $a_e'$ be the two endpoints of the arcs $m_e$ and $b_e$ and $b_e'$ be the endpoints of the arcs $n_e$.  The points $a_e, a_e', b_e, b_e'$ induce points on the boundary of $G$, $\widehat{G}$, and $H_e$, and on the states $\hat{s}, s$, and $s_e$. Now consider a boundary component of $s$. If this boundary component does not contain any of the points $a_e, a_e', b_e, b_e'$ then there is a naturally corresponding boundary component in $\hs$. If the boundary cycle contains any of the points $a_e, a_e', b_e, b_e'$ then there is a corresponding boundary cycle in $\hs$ that contains the same set of points. This sets up a natural correspondence between the boundary components in $s$ and a set of boundary components of $\hs$. 
However, $\hs$ has extra boundary components arising from the states
$s_e \in \cH \cup \dH$. These extra boundary components  are
precisely the unmarked  boundary components  in the states $s_e \in
\cH \cup \dH$, $e\in E$.
Since the    points $a_e, a_e', b_e, b_e'$ belong to the same boundary
component if and only if they belong to the same connected component,
we have that for each $e\in E$ for which $s_e \in \cH$, there are
$\partial (s_e) -1$ unmarked boundary components; and for 
each $e\in E$ for which $s_e \in \dH$, there are   $\partial (s_e)
-2$ unmarked boundary components.
Therefore 
\[  \partial (\hs) = \partial (s) +\sum_{\substack{e\in E \\ s_e \in
\cH}}  (\partial (s_e) -1) 
+\sum_{\substack{e\in E\\ s_e \in \dH}}  (\partial (s_e) -2)  . \]
The lemma then follows.
\end{proof}

\subsection{An expansion for the \BR polynomial} \label{sec:planBR}
We will consider the following state sums:
\begin{equation}\label{eq:planOGS}
\begin{array}{l}
\ceta (a,b,c)  := \sum_{s \in \cH}
a^{k(s)-1} b^{e(s)}c^{\partial (s)-1}, \\
\deta (a,b,c)  := \sum_{s \in \dH}
a^{k(s)-2} b^{e(s)}c^{\partial (s)-2}. 
\end{array}\end{equation}

The following lemma is analogous to Lemma~\ref{lem:tutteOGS}.
\begin{lemma}
Let $\left(G , \{ H_e \}_{e\in E}\right)$ be a 2-decomposition of
$\hG$ then
\begin{equation}\label{eq:BRogs}
Z(\hG ; a,b,c) =
\left( \prod_{e\in E}( \deta) \right)
Z\left(G; a, \left\{ \ceta / \deta  \right\}_{e \in E}
,c  \right).
\end{equation}
\end{lemma}
The proof of this lemma is a direct generalization of the proof of
Lemma~\ref{lem:tutteOGS} (using the additional relation
$\partial(\hs) =  \partial(s)+\sum (\partial(s_e)-2) + \sum (\partial(t_e)-1)$ 
arising from Lemma~\ref{lem:brcon}), and is therefore omitted.

%
%
%

To find a formula for $Z(\hG)$ it remains to determine $\deta$ and
$\ceta$.

Since the edge set of $H_e$ is a subset of the edge set of $A_e$,
we may view states of $H_e$ as states of $A_e$. The states of
$A_e$ can then be partitioned into four subsets
\begin{equation}\label{eq:planBRpartition}
\begin{array}{ll}
 \cH, \quad &  \mathcal{T}^1 (H_e) := \{ s \cup e |  s \in
\cH \} \\
\dH,
& \mathcal{T}^2 (H_e):= \{ s \cup e |  s \in \dH \} .
\end{array}\end{equation}
Consider the effect of the insertion of the edge $e$ in a state
of $H_e$ on the numbers of connected components, edges, and boundary
cycles. There are two cases. If $s_e \in \cH$ then the insertion of
$e$  increases the number of boundary cycles by one, so
that if $s_e$ contributes  the term   $a^k b^e c^{\partial}$ to the
\BR polynomial then the state obtained by inserting $e$ contributes 
$ a^k (b^e x_e) c^{\partial +1} $. If $s_e \in \dH$ then the insertion
of $e$ also decreases the number of boundary cycles by one and the
number of connected components by one. This means that if $s_e$
contributes  the term   $a^k b^e c^{\partial} $ to the \BR polynomial
then the state obtained by inserting $e$ contributes  $ a^{k-1} (b^e
x_e) c^{\partial-1} $.

We can now separate the terms in $Z(A_e; a, \vb ,c)$  (with $\vb$ as
in Notation~\ref{notation1}) arising from the four subsets in the
partition and write
\begin{equation*}
Z(A_e; a, \vb ,c) 
=
(x_e c + 1)
 \sum_{s \in \cH}a^{k(s)}b^{e(s)}c^{\partial(s)} 
+ 
(x_e a^{-1}c^{-1} + 1)
 \sum_{s \in \dH} a^{k(s)}b^{e(s)}c^{\partial(s)}  ,
\end{equation*} 
or
\[
Z(A_e; a, \vb ,c)  =
(x_e c + 1) ac\; \ceta
 +
(x_e a^{-1}c^{-1} + 1) a^{2}c^{2} \; \deta.
\]
Writing this as a linear equation in $x_e$, and using the
deletion-contraction relation (\ref{eq:delcont}), we have
\begin{equation}\label{eq:plan}
ac   \left( \ceta + ac\deta \right)
+
ac
\left(  c\ceta +  \deta \right)x_e
=
Z(H_{e};a,b,c)+x_{e}Z(A_{e}/e;a,b,c).
\end{equation}
This gives  rise to a system of equations
\begin{equation}\label{eq:plansolv}
\begin{split}
ac   \left(  \ceta +ac\deta  \right) &= Z(H_{e};a,b,c) \\
ac
\left(  c \ceta +  \deta  \right)
&=
Z(A_{e}/e;a,b,c).
\end{split}
\end{equation}
This pair of linear equations uniquely determines  $\ceta$
and $\deta$. Substitution into equation~(\ref{eq:BRogs}) then gives
the following theorem.
\begin{theorem}
Let $\hG$ be an embedded graph with a 2-decomposition $\left( G , \{
H_e \}_{e\in E} \right)$, such that each graph $H_e$ is embedded in a
neighbourhood of the edge $e$ of the embedded graph $G$. In addition
let $A_e$ be the ribbon graph $H_e$ with an additional ribbon $e$
joining the vertices $u_e$ and $w_e$.
 Then
\[ Z (\hG;a, b ,c) =
(ac)^{-e(G)}
\left( \prod_{e\in E} g_e  \right)
Z\left(G; a, \left\{  f_e / g_e \right\}_{e\in E},  c
\right).
\]
where  $f_e$ and $g_e$ are the solutions to
\begin{align*}
acg_e+f_e&=Z(H_{e};a,b,c)\\
g_e+cf_e&=Z(A_{e}/e;a,b,c).
\end{align*}
\hfill$\Box$
\end{theorem}

We may use this to find a result analogous to
Corollary~\ref{cor:bry}.
\begin{corollary}\label{cor:planbryII}
Let $G=(V,E)$ be a ribbon graph, $A$ be a planar ribbon graph, and
$H=A\backslash e$.
Then
\[
R(G \otimes A; \al, \be  ,\ga) =
(h)^{n(G)} (h^{\prime})^{r(G)}
R\left(G;  \frac{R( H ; \al,\be ,\ga )}{h^{\prime} }, \frac{\be
h^{\prime} }{ h} ,
\ga \right),
\]
where $h$ and $h^{\prime}$ are the unique solutions to
\begin{align*}
 h + \be h^{\prime} &=  R(A/e; \al, \be,\ga ) \\
(\al-1)h + h^{\prime} &= R( H ; \al, \be,\ga ).
\end{align*}
\end{corollary}
\begin{proof}
The proof is similar to that of Corollary~\ref{cor:bry}. By
(\ref{eq:RZ}) we have 
 \[
 R(\goa ; \al, \be, \ga) = (\al -1)^{-k(\goa)} (\be\ga)^{-v(\goa)}
Z\left(\goa ; (\al-1) \be \ga^2 ,\be \ga, \ga^{-1} \right).
 \] 
 An application of  the theorem gives
 \[
 R(\goa ; \al, \be, \ga)  = (\al -1)^{-k(\goa)} (\be\ga)^{-v(\goa)}
g^{e(G)}
 Z\left(G ; (\al-1) \be \ga^2 ,f/g, \ga^{-1} \right),
 \] 
  where $f$ and $g$ are the solutions to   
  \begin{align*}
  (\al-1)\be\ga ( f+ (\al-1)\be\ga g  ) &= Z\left(H; (\al-1) \be
\ga^2 ,\be \ga, \ga^{-1} \right),
  \\
    (\al-1)\be\ga (\ga^{-1} f+  g  )  &= Z\left(A/e; (\al-1) \be
\ga^2 ,\be \ga, \ga^{-1} \right).
  \end{align*}
(Note that the $(ac)^{-e(G)}=((\alpha-1)\beta\gamma)^{-e(G)}$ factor
has been incorporated into these equations.)
  
By a second application of (\ref{eq:RZ}) we can rewrite the above as
\[ R(\goa ; \al, \be, \ga) 
=
(\al -1)^{-k(\goa)} (\be\ga)^{-v(\goa)} g^{e(G)}
\left( \frac{(\al-1) \be \ga g }{f}\right)^{k(G)} \left(\frac{f}{g}
\right)^{v(G)} 
R\left(  G; \frac{(\al-1) \be \ga g +f}{f}  , \frac{f}{\ga g}  , \ga
\right).
\]
Then, just as in the proof of Corollary~\ref{cor:bry}, making the
substitution
\[
h:=(\al-1)^{-k(H)+1}(\be\ga)^{-v(H)+2} g, \quad h^{\prime}:=
(\al-1)^{-k(H)+1}(\be\ga)^{-v(H)+1} f,
\] 
remembering that $v(\goa)=(v(H)-2)e(G)+v(G)$ and
$k(\goa)=(k(H)-1)e(G)+k(G) $, and cancelling terms, we have
 \[
 R(G \otimes A; \al, \be ,\ga) =
(h)^{n(G)} (h^{\prime})^{r(G)}
R\left(G;  ((\al-1)h + h^{\prime} )/h^{\prime} , \be h^{\prime} / h ,
\ga \right).
 \]
Now using (\ref{eq:RZ}) and noting that $k(H)=k(A/e)$ and
$v(A/e)=v(H)-1$, the linear equations can be written as 
\[
 h + \be h^{\prime} =  R(A/e; \al, \be,\ga ), \quad
(\al-1)h + h^{\prime} = R( H ; \al, \be,\ga ).
\]
\end{proof}

We use this corollary to extend our Example~\ref{ex:1}.
\begin{example}
Let $A = C_3$, the 3-cycle. Then $h=\al +1$ and $h^{\prime} = 1$ and
we have
$R(G \otimes C_3;\al,\be,\ga)=(\al +1)^{n(G)} R(G;  \alpha^2, \be /(\al+1) ,\ga   )$. See
also \cite{Mo} for a proof of this fact using knot theory.
\end{example}

The following is an important application of our results. It provides
a method for constructing infinitely many pairs of distinct ribbon
graphs with the same \BR polynomial.
\begin{corollary}
Let $G\otimes H$ and $G^{\prime}\otimes H$ be two embedded graphs
with the property that each copy of $H$ is embedded in the
neighbourhood of an edge. Then if 
 $R(G; a,b ,c )= R(G^{\prime};a,b,c)$,  $R(G\otimes H; a,b ,c )=
R(G^{\prime}\otimes
H;a,b,c)$.
\end{corollary}

\begin{proof}
This follows from the previous corollary because if
$R(G)=R(G^{\prime})$ then, by setting $c=1$, $T(G)=T(G^{\prime})$, and
since the rank and nullity of a graph can be recovered from its Tutte
polynomial, we have $ r(G)=r(G^{\prime})$ and $n(G)=n(G^{\prime})$.
\end{proof}


\section{The \BR polynomial II: the general case}\label{sec:brII}

We begin with an informal discussion of the main ideas in this
section.
Consider the construction of a ribbon graph $\hG$ from the
2-decomposition $\left( G , \{ H_e \}_{e\in E} \right)$ locally at an
edge $e=(u,w)$ of the template $G$. We will think of the construction
of $\hG$ as the identification of the marked points $m$ and $n$ on
the  vertices $u$ and $w$ of $H_e$ with the corresponding marked
points
$m$ and $n$ on the vertices $u$ and $w$ of the template $G\bs e$. 

Begin by partitioning  $\s(H_e)$ according to the boundary and
connected components containing the marked points $m$ and $n$ on
the vertices $u$ and $v$: we let $\ddH$ be the set of states of $H_e$
in which $m$ and $n$ lie in different connected and different boundary
components; $\bcH$ the states in which $m$ and $n$ lie in the same
connected and same  boundary components; and $\dcH$ the states in
which
$m$ and $n$ lie in the same connected component but different boundary
components.

We would like to define a replacement operation on $\s(G)$ to
construct $\s (\hG)$, but we run into a problem. The edge $e$ is
either in a state of $G$, in which case  we can glue in states from
$\bcH$ (reflecting the fact that $m$ and $n$ lie in the same connected and same boundary components), or $e$ is not in a state of
$G$, in which case  we can glue in states from $\ddH$ (reflecting  the fact that $m$ and $n$ lie in different connected and different boundary components). No states of $G$ ever reflect the fact
that the markings $m$ and $n$ in $H_e$ can lie in different boundary
components  and the same connected component, so in this construction
we never glue in any states from $\dcH$.

To get around this problem  we replace each edge $e$ in $G$ with an
edge $g_e$ and a loop $f_e$ as in Figure~\ref{fig:gtilde}, to obtain
a graph
$\tG$. Then there is a subset of the edges $\{f_e, g_e\}$ for each
choice of $e\in E$, such that the connectivity and boundary
connectivity of the markings $m_e$ and $n_e$ in $\tG$ correspond to
the  connectivity and boundary connectivity of the markings $m_e$ and
$n_e$ in each of $\ddH$, $\dcH$ and $\bcH$. Define a replacement
operation as follows: $\{g_e\}$ is replaced by states from $\bcH$;
$\{f_e , g_e\}$ by states from $\dcH$; and the states $\emptyset$ and
$\{f_e\}$ are replaced by states from $\ddH$. The considerations
above show that every element of $\s (\hG)$ is obtained from
$\s(\tG)$ by this operation. However, the states are  not obtained
uniquely: we have two configurations,  $\emptyset$ and $\{f_e\}$, for
which we substitute states from $\ddH$. So rather than dealing with
$\s (\tG)$ we will deal with equivalence classes in this set in which
we identify the states containing $f_e$ but not $g_e$, or neither
$f_e$ or $g_e$. Our replacement operation will then
give a construction of $\s (\hG)$  from the sets $\s (\tG)/\sim$ and
$\ddH$, $\dcH$ and $\bcH$, $e\in E$.
We will now do this formally.

\subsection{Decomposing the ribbon graph} \label{ss:brCyDe}

We view the construction of $\hG$ in terms of the identification of
arcs $n_e$ and $m_e$ as described in Subsection~\ref{ss:CyDe} and by
Figure~\ref{fig:dmark}.

\begin{figure}
\[\epsfig{file=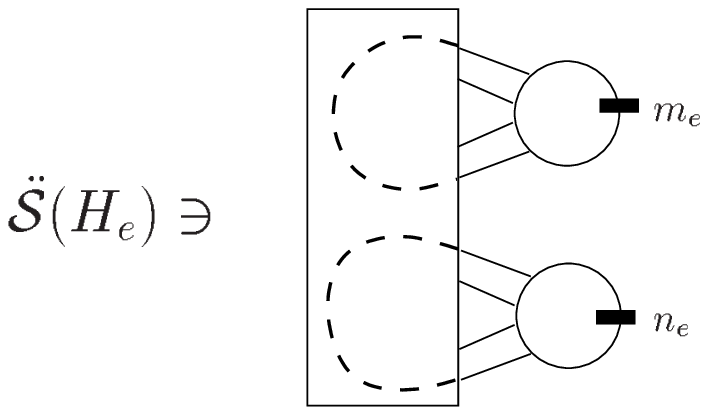, width=4cm} \;\; , \hspace{2cm}
\epsfig{file=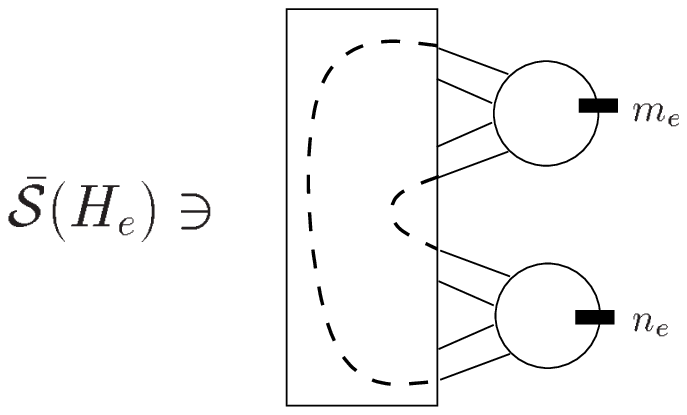, width=4cm}\]
\caption{}
\label{fig:CyDe2}
\end{figure}

Given a state $\hs \in \s (\hG)$ and a 2-decomposition $\left(G, \{
H_e
\}_{e\in E} \right)$ for the ribbon graph  $\hG$, the state $\hs$
uniquely determines a set of states $\{ s_e \in \s (H_e)\}_{e\in E}$.
Partition each $\s (H_e)$ as follows. Let $s_e \in \s(H_e)$, then 
$s_e \in \bar{\s}(H_e)$ if and only if $m_e$ and $n_e$ belong to the
same boundary component. Otherwise $s_e \in \ddot{\s}(H_e)$. (So the
accent
on $\s$ has one component if and only if $m_e$ and $n_e$  belong to
one boundary component.) These two situations are indicated in
Figure~\ref{fig:CyDe2}.

We want to include connectivity information in the above partition.
We do this by taking its intersection with the partition of
Subsection~\ref{ss:TuDe}.
Let $\left( G , \{ H_e \}_{e\in E} \right)$ be a 2-decomposition
for the ribbon graph $\hG$ and
$\cH$, $\dH$, $\ddot{\s}(H_e)$ and $\bar{\s}(H_e)$ be the sets
described in Subsections~\ref{ss:TuDe} and above. Define
\[\begin{array}{ll}
 \bar{\s}^1 (H_e)  :=  \cH \cap \bar{\s}(H_e), \quad &
 \ddot{\s}^1(H_e):=  \cH \cap \ddot{\s}(H_e), \\ 
\bar{\s}^2(H_e):=   \dH \cap \bar{\s}(H_e), \quad &
\ddot{\s}^2(H_e) :=  \dH \cap \ddot{\s}(H_e).
\end{array}\]
Notice that $\bdH = \emptyset$ and
$\dcH \cup \ddH = \ddot{\s}(H_e)$.
We have
\begin{lemma}\label{lem:BRdecomp}
$\bcH$, $\dcH$ and $\ddH$
partition the sets $\s (H_e)$, $e \in E$.
Moreover every state of $\hG$
can be uniquely obtained by the replacement of an edge in a state of
$G$ by an element of $\bcH$ and the replacement of a
non-edge in that state of $G$ with an element of $\dcH
\cup \ddH$.
\hfill$\Box$
\end{lemma}

\medskip

\begin{figure}
\begin{center}
\subfigure[]{\epsfig{file=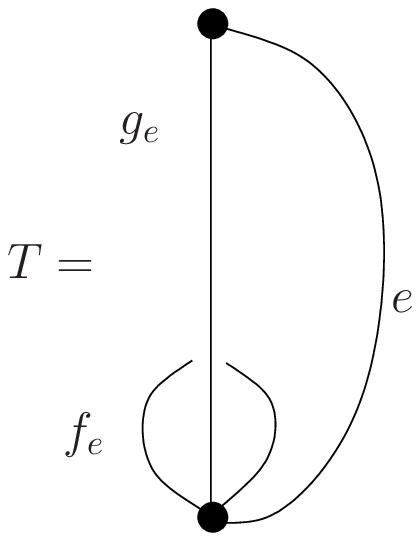, height=3cm}} \quad \quad \quad
\subfigure[]{\epsfig{file=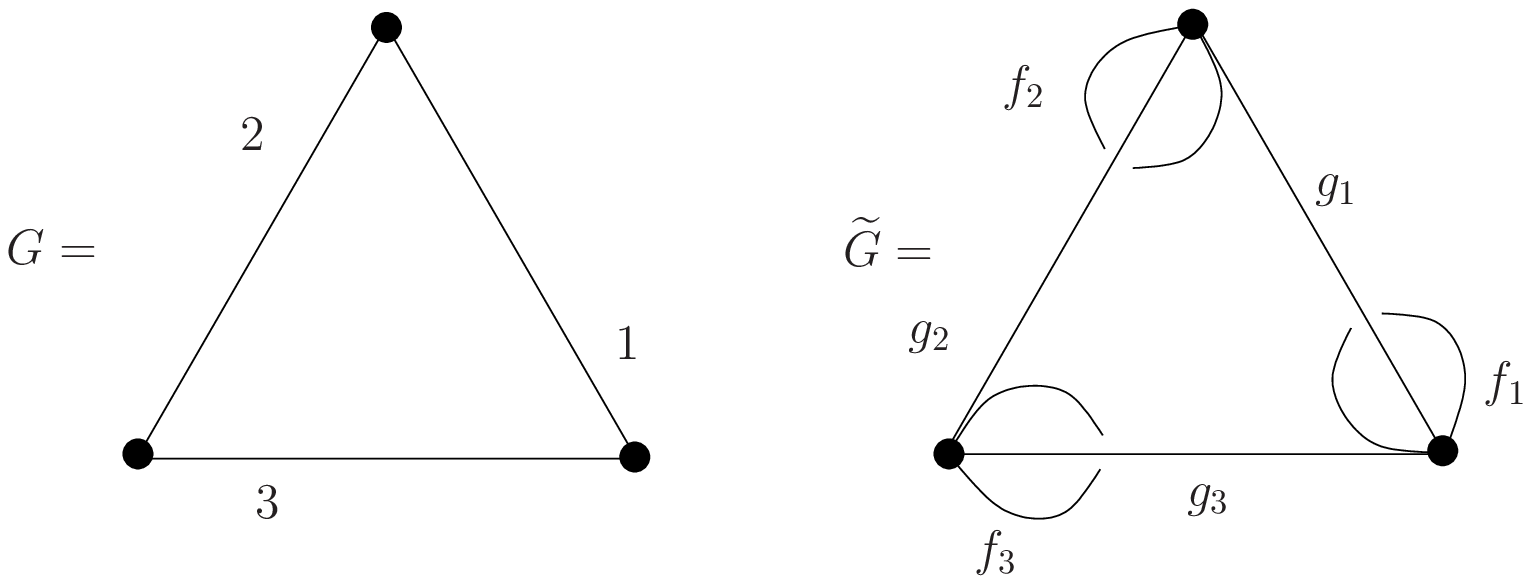, height=3cm}}
\end{center}
\caption{}
\label{fig:gtilde}
\end{figure}

Rather than considering states of the template $G$, we need
to consider states of a slightly more complex ribbon graph. We define
the ribbon graph $\tG$ to be the tensor product (along $e$) of the
template $G$ with the ribbon graph $T$ defined by
Figure~\ref{fig:gtilde}(a),
having the specified edge-labels. We will see that it does not matter
if
the loop is at the vertex labelled $u_e$ or $w_e$.   
An example of $\tG$ is shown in Figure~\ref{fig:gtilde}(b).

We want to construct the set of states $\s (\hG)$ by replacing states 
in $\s (\tG)$ by states from $\s(H_e)= \bcH \cup \dcH\cup \ddH$.


We will say that two states of $\tG = G \otimes T$ are
equivalent if for
some choice of  $e$, one state contains neither of the edges $f_e$ or
$g_e$,
the other state contains the edge $f_e$ but not $g_e$, and  all of
the other
edges contained in the two states are the same. This defines an
equivalence
relation $\sim$ on $\s (\tG)$. We are interested in the set $
\s(\tG)/\sim$.

We can construct the set of states $\s (\hG)$ by replacing
equivalence classes  of
$ \s (\tG)/\sim$ with $\{ \s (H_e) \}_{e \in E(G)}$ as follows. Given
an
equivalence class $[s] \in \s (\tG)/\sim$ choose a representative $s
\in [s]$.
Then for each choice of $e$, we have the following possibilities in
the equivalence
class: $s$ contains both $f_e$ and $g_e$; or $g_e$ but not $f_e$; or
it does
not contain $g_e$ and may or may not contain $f_e$ (these two
situations 
being equivalent).

\begin{figure}
\[\epsfig{file=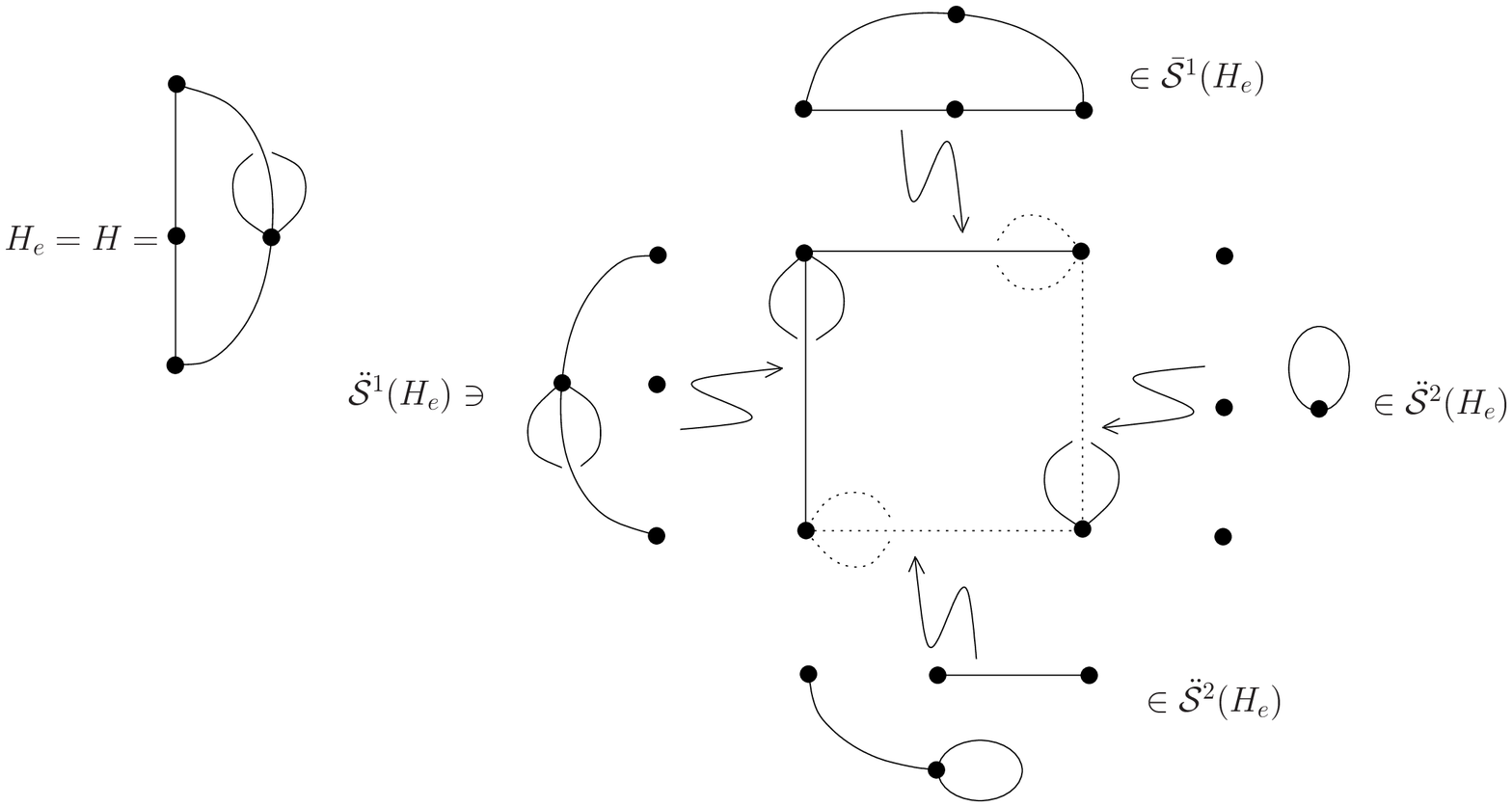, width=15cm}\]
\caption{}
\label{fig:Cydecomp}
\end{figure}

Construct a state of $\hG$ using the  replacement operation: 
 \begin{itemize}
 \item if $s$ contains both of $f_e$ and $g_e$, remove both of these
two edges and glue in a state from $\dcH$;
 \item if a state contains the edge $g_e$ but not $f_e$, remove the
edge $g_e$ and glue in a state from  $\bcH$;
 \item if a state contains the edge $f_e$ but not $g_e$, remove
the edge $f_e$ and glue in a state from $\ddH$, or if the original
state
contains neither of the edges $f_e$ or $g_e$ then glue in a state from
$\ddH$.
 \end{itemize}
An example  of the decomposition of a state of $\hG$ is shown in
Figure~\ref{fig:Cydecomp} where $G=C_4$.

The following lemma is clear.
\begin{lemma}\label{lem:BRdecompII}
$\bcH$, $\dcH$ and $\ddH$
partition the sets $\s (H_e)$, $e \in E$.
Moreover every state of $\hG$
can be uniquely obtained by replacing classes in $\s(\tG)/\sim$ with 
elements of $\bcH$, $\dcH$ and $\ddH$ in the manner described above.
\hfill$\Box$
\end{lemma}

\subsection{An expansion for the \BR polynomial} \label{sec:BR}

We consider the following state sums:
\begin{equation}\label{eq:ribOGS}
\begin{array}{l}
\Phi_{\tG}(a, \{f_e,g_e\}_{e\in E} ,c)  := 
\sum_{[s] \in \s(\tG)/\sim} a^{k(s)} \left( \prod_{e \in s}x_e
\right)   c^{\partial (s)}, \\
\doeta (a,b,c) :=  \sum_{s \in \dcH}
a^{k(s)-1} b^{e(s)} c^{\partial (s)-2}, \\
\dteta (a,b,c)  := \sum_{s \in \ddH}
a^{k(s)-2} b^{e(s)}c^{\partial (s)-2}, \\
\boeta (a,b,c)  := \sum_{s \in \bcH}
a^{k(s)-1} b^{e(s)}c^{\partial (s)-1},
\end{array}\end{equation}
where in the expression for $\Phi_{\tG}$ the product is taken over the
edges of a representative $s\in [s]$ such that $s$ has the fewest
edges
in its equivalence class, and $\{x_e\}$ denotes the set of labels
$\{f_e ,g_e\}$ of the edges in $\tG$.  

We emphasize that we are abusing notation by using $\{f_e ,g_e\}$ to denote both the set of edge labels and the set of edge weights. This abuse of notation should cause no confusion.

\begin{lemma}\label{l.referee}
Suppose that a state $\hs$ of $\hG$ is obtained by replacement
from the states
$[s] \in \s(\tG)$, $s_e \in \dcH$, $t_e \in \ddH$, and $u_e \in\bcH$ 
using the decomposition. Then
\begin{equation*}
\begin{array}{ll}
\partial(\hs)  = & \partial (s) + \sum (\partial(s_e)-2) + \sum
(\partial(t_e)-2)+
\sum (\partial(u_e)-1),
\\
k(\hs)  =& k(s) + \sum (k(s_e)-1) + \sum (k(t_e)-2) + \sum (k(u_e)-1),
\\
e(\hs) =&  \sum e(s_e) + \sum e(t_e)+ \sum e(u_e),
\end{array}
\end{equation*}
where the representative $s\in [s]$ is chosen so that it has the
fewest edges in its class.
\end{lemma}
\begin{proof}
For the first identity, since our representative $s\in [s]\in
\s(\tG)/ \sim$ was chosen so that it has the fewest edges in its
class,
the state $s$ contains an edge $f_e$ if and only if it contains the edge $g_e$. The argument now
follows
the proof of Lemma~\ref{lem:brcon}: 
each boundary component of $s \in \s
(\tG)$ corresponds to a boundary component of $\hs \in \s (\hG)$.
These
are easily seen to be the boundary components of $\hs$ which contain a
point $a_e, a_e', b_e, b_e'$ (with these points defined as in Lemma~\ref{lem:brcon}). 
However $\hs$ has additional boundary components arising from the
states $s_e$. These extra boundary components of $\hs$ are precisely
the boundary components of the $s_e$ which do not contain a
point $a_e, a_e', b_e, b_e'$. For each $e$ there are $(\partial(s_e)-2)$ of
these
if $s_e \in \dcH$, 
$(\partial(t_e)-2))$ of these if $t_e \in  \ddH$, and
$(\partial(u_e)-1))$ of these if $u_e \in  \bcH$. The result follows. 

The remaining identities follow similarly (see also the proofs of
Lemmas \ref{lem:con} and \ref{lem:brcon}).
\end{proof}

We define the linear map
$\mathcal{F}: \mathbb{Z} [\{f_e, g_e\}_{ e\in E}] \rightarrow
\mathbb{Z}[\{\doeta, \dteta, \boeta\}_{e\in E}]$
 to
be the linear extension of the map
\[
\mathcal{F} : \prod_{e\in E} f_e^{\alpha_e}g_e^{\beta_e} \mapsto
\prod_{e\in E} \left( \doeta  \right)^{\alpha_e \beta_e}
\left( \dteta \right)^{(1-\beta_e)}
\left( \boeta \right)^{(\beta_e-\alpha_e \beta_e)} .
\]
For example, $\mathcal{F}$ sends the monomial $f_1g_1g_2 \in
\mathbb{Z}[f_1,f_2,f_3,g_1,g_2,g_3 ]$ to 
$\ddot{\eta}_1^{(1)}\bar{\eta}_2^{(1)}\ddot{\eta}_3^{(2)}$.

\begin{lemma}\label{lem:zfphi}
\[
Z( \hG ; a,b,c) = \mathcal{F} \left(
\Phi_{\tG} 
 \right) . 
\]
\end{lemma}
\begin{proof}
Let $[s]\in \s (\tG)/\sim$ and $s\in [s]$ be the representative with
the fewest edges in its class. We then know that for each index $e\in
E$, $s$ contains both $f_e$ and $g_e$; or $g_e$ but not $f_e$; or
neither $f_e$ nor $g_e$. We construct a state of $\hG$ by replacing
the three edge configurations $\{f_e,g_e\}$, $\{g_e\}$, and
$\emptyset$
of the pair of edges $\{f_e, g_e\}$, with states $s_e$ of $\dcH$,
$\bcH$, and $\ddH$ respectively. Each state of $\hG$ is uniquely
obtained in this way, and the corresponding contribution to $Z(\hG)$
will
be
  \begin{equation}\label{eq:co} 
  a^{k(s)}c^{\partial(s)}\prod_{e\in E}
  \left( a^{k(s_e)-1}b^{e(s_e)}c^{\partial (s_e)-2}
\right)^{\alpha_e\beta_e}
    \left( a^{k(s_e)-2}b^{e(s_e)}c^{\partial (s_e)-2}
\right)^{1-\beta_e}
      \left( a^{k(s_e)-1}b^{e(s_e)}c^{\partial (s_e)-1}
\right)^{\beta_e-\alpha_e\beta_e}
  ,\end{equation}
 where 
$  \alpha_e = \left\{ \begin{array}{ll} 1 & \text{ if } f_e \in s \\
0 & \text{ otherwise}  \end{array}\right.  $
and 
$  \beta_e = \left\{ \begin{array}{ll} 1 & \text{ if } g_e \in s \\ 0
& \text{ otherwise}  \end{array}\right.  $.
 
Clearly expression~(\ref{eq:co}) is equal to $\mathcal{F} \left(
a^{k(s)} (\prod_{e\in s} x_e ) c^{\partial (s)})
\right)$, the contribution of the state $s$ to  $\Phi_{\tG}$. But it
is also equal to
 \begin{multline*}
 a^{ k(s) + \sum_{s_e \in \dcH} (k(s_e)-1)+\sum_{s_e \in \ddH}
(k(s_e)-2)+\sum_{s_e \in \bcH} (k(s_e)-1)} \\
 b^{  \sum_{s_e \in \dcH} e(s_e)+\sum_{s_e \in \ddH} e(s_e)+\sum_{s_e
\in \bcH} e(s_e)} \\
 c^{ \partial(s) + \sum_{s_e \in \dcH} (\partial(s_e)-2)+\sum_{s_e
\in \ddH} (\partial(s_e)-2)+\sum_{s_e \in \bcH} (\partial(s_e)-1)}. 
   \end{multline*}
From Lemma~\ref{l.referee}, this sum is equal to
$a^{k(\hs)}b^{e(\hs)}c^{\partial (\hs)}$. 
By the uniqueness of the decomposition of the states $\hs \in \s
(\hG)$ into states of $\s(\tG)/\sim$ and $\s (H_e)$, $e\in E$, the
result follows on summing over the states.
\end{proof}

Notice that the polynomial $Z(\tG; a, \vx , c)$ enumerates all of the
states of $\tG$. Our next result uses this observation to replace $
\Phi_{\tG}$ in the above lemma with $Z(\tG)$.

\begin{lemma}\label{lem:zgz}
Let $\left(G , \{ H_e \}_{e\in E}\right)$ be a 2-decomposition of
$\hG$, and 
$\mathcal{G}: \mathbb{Z} [\{f_e, g_e\}_{e\in E}] \rightarrow
\mathbb{Z}[\{\doeta, \dteta, \boeta\}_{e\in E}]$  be the linear
extension
of the map
\[
\mathcal{G} : \prod_{e\in E} f_e^{\alpha_e}g_e^{\beta_e} \mapsto
\prod_{e\in E} \left( \doeta  \right)^{\alpha_e \beta_e}
\left( \frac{1}{2}\dteta \right)^{(1-\beta_e)}
c^{(\alpha_e \beta_e-\alpha_e)}
\left( \boeta \right)^{(\beta_e-\alpha_e \beta_e)}.
\]
Then
\[  Z(\hG ; a,b,c) = \mathcal{G}\left(  Z(\tG; a, \vx , c) \right).\]
\end{lemma}
\begin{proof}
We can write the polynomial 
\[  Z(\tG; a, \vx , c)    =    \sum_{\ts \in \tG} a^{k(\ts)}
c^{\partial (\ts)} \left( \prod_{e\in \ts } x_e\right)\]
as
\[   \sum_{\ts \in \tG} a^{k(\ts)} c^{\partial (\ts)}  \prod_{e\in
\ts } \left(  (f_eg_e)^{\alpha_e}
(f_e)^{\beta_e}(g_e)^{\gamma_e}(1)^{\delta_e}\right), \]
where exactly one of $\alpha_e$, $\beta_e$, $\gamma_e$, $\delta_e$
is one and all of the others are zero for each $e\in E$.

Now suppose that a state $\ts$ decomposes into states $s_e \in \s
(H_e)$, $e\in E$. Then for any $e$, if $\beta_e=1$ or $\delta_e=1$ we
know that $s_e \in \ddH$ and since a state containing only $f_e$ will
have exactly one more boundary component than a state containing
neither $f_e$ or $g_e$, for some $e$, we may write the above formula
as
\begin{equation}\label{eq:gproof1} \sum_{[\ts] \in \tG/\sim}
a^{k(\ts)} c^{\partial (\ts)}
\prod_{e\in \ts } \left(  (f_eg_e)^{\alpha_e}
(1+cf_e)^{\beta_e}(g_e)^{\gamma_e}\right), \end{equation}
where the representative $\ts$ is chosen so that it has the fewest
edges in its class. 

We need to show that the map $\mathcal{G}$ applied to
(\ref{eq:gproof1}) is equal to $\mathcal{F}(\Phi_{\tG})$. The state
sum $\Phi_{\tG}$ can be expressed as
\begin{equation}\label{eq:gproof2}
\Phi_{\tG} = \sum_{[\ts ] \in \tG/\sim} a^{k(\ts)} c^{\partial (\ts)}
\prod_{e\in \ts } \left(  (f_eg_e)^{\alpha_e}
(1)^{\beta_e}(g_e)^{\gamma_e}\right),
\end{equation}
 where the representative $\ts$ is chosen so that it has the fewest
edges in its class, and exactly one of $\alpha_e$, $\beta_e$,
$\gamma_e$, $\delta_e$
is one and all of the others are zero for each $e\in E$.

There is a clear correspondence between the summands  of
(\ref{eq:gproof1}) and (\ref{eq:gproof2}). 
In particular, if for some $e$, a summand of (\ref{eq:gproof1})
contains the expression $(f_eg_e)^{1}(1+cf_e)^{0}(g_e)^{0}$  then the
corresponding summand of (\ref{eq:gproof2}) contains a term
$(f_eg_e)^{1}(1)^{0}(g_e)^{0}$ and these terms are mapped by
$\mathcal{G}$ and $\mathcal{F}$  respectively to $\doeta$.
Also if  for some $e$, a summand of (\ref{eq:gproof1}) contains the
expression $(f_eg_e)^{0}(1+cf_e)^{0}(g_e)^{1}$  then the
corresponding summand of (\ref{eq:gproof2}) contains a term
$(f_eg_e)^{0}(1)^{0}(g_e)^{1}$ and these terms are mapped by
$\mathcal{G}$ and $\mathcal{F}$  respectively to $\boeta$.
Finally, if for some $e$, a summand of (\ref{eq:gproof1}) contains
the expression $(f_eg_e)^{0}(1+cf_e)^{1}(g_e)^{0}$  then the
corresponding summand of (\ref{eq:gproof2}) contains a term
$(f_eg_e)^{0}(1)^{1}(g_e)^{0}$. In this case 
\[ \mathcal{G}((f_eg_e)^{0}(1+cf_e)^{1}(g_e)^{0}) 
= \frac{1}{2} \dteta +c(c^{-1}\frac{1}{2} \dteta  ) = \dteta =
\mathcal{F}((f_eg_e)^{0}(1)^{1}(g_e)^{0}).\]
Hence we see that applying the map $\mathcal{G}$
to (\ref{eq:gproof1}) 
will give $\mathcal{F}(\Phi_{\tG})$, and then an application of
Lemma~\ref{lem:zfphi} will give the required identity
\[    \mathcal{G}\left(  Z(\tG; a, \vx , c) \right) =Z(\hG ; a,b,c)
.\] 
\end{proof}

\medskip

It remains to determine $\boeta$ , $\doeta$ and  $\dteta$.

\subsection{Using ribbon graphs}

We may view states of $H_e$ as states of $A_e$. The states of
$A_e$ can be partitioned into six subsets
\begin{equation}\label{eq:BRpartition}
\begin{array}{ll}
 \dcH, \quad & \ddot{\mathcal{T}}^1 (H_e) := \{ s \cup e |  s \in
\dcH \} \\
\ddH,
& \ddot{\mathcal{T}}^2(H_e) := \{ s \cup e |  s \in \ddH \} \\
\bcH, & \bar{\mathcal{T}}^1(H_e) :=\{ s \cup e |  s \in
\bcH \}.
\end{array}\end{equation}

Consider the effect of the insertion of the edge $e$ into a state
of $H_e$ on the number of connected components, edges, and boundary
components, and the corresponding terms in $Z(H_e;a,\vb,c)$. There
are three cases. If $s_e \in \dcH$ then the insertion of $e$ 
decreases the number of boundary cycles by one. This means
that if $s_e$ contributes the term $a^k b^e c^{\partial}$ to the
\BR polynomial then the state obtained by inserting $e$ contributes
$ a^k (b^e x_e) c^{\partial-1} $. If $s_e \in \ddH$ then the insertion
of $e$  decreases the number of boundary cycles by one and the
number of connected components by one. This means that if $s_e$
contributes  the term   $a^k b^e c^{\partial} $ to the \BR polynomial
then the state obtained by inserting $e$ contributes $ a^{k-1} (b^e
x_e) c^{\partial-1} $. Finally, if $s_e \in \bcH$ then the insertion
of
$e$ also increases the number of boundary cycles by one. This means
that if $s_e$ contributes  the term $a^k b^e c^{\partial} $ to the
\BR polynomial then the state obtained by inserting $e$ contributes
$ a^{k} (b^e x_e) c^{\partial +1} $.

We can then separate the terms in $Z(A_e; a, \vb ,c)$, $\vb$ as in
Notation~\ref{notation1}, arising from the six subsets in the
partition and write:
\begin{multline*}
Z(A_e; a, \vb ,c) 
=
(x_e c^{-1} + 1) \sum_{s \in \dcH}
a^{k(s)}b^{e(s)}c^{\partial(s)} + 
(x_e a^{-1}c^{-1} + 1) \sum_{s \in \ddH}
a^{k(s)}b^{e(s)}c^{\partial(s)}  \\  + 
(x_e c + 1) \sum_{s \in \bcH}
a^{k(s)}b^{e(s)}c^{\partial(s)}.
\end{multline*} 
Rewriting in terms of $\doeta,\dteta$, and $\boeta$ gives
\[
Z(A_e; a, \vb ,c)  =
(x_e c^{-1} + 1) ac^{2}\; \doeta
 +
(x_e a^{-1}c^{-1} + 1) a^{2}c^{2} \; \dteta +
(x_e c + 1)ac \; \boeta.
\]
Writing this as a linear equation in $x_e$, and using the
deletion-contraction relation (\ref{eq:delcont}), we have
\begin{equation}\label{eq:fullunsolv}
ac   \left( c \doeta +ac\dteta + \boeta \right)
+
ac
\left(  \doeta +  \dteta + c \boeta \right)x_e
=
Z(H_{e};a,b,c)+x_{e}Z(A_{e}/e;a,b,c).
\end{equation}
giving rise to a system of equations
\begin{equation}\label{eq:unsolv}
\begin{split}
ac   \left( c \doeta +ac\dteta + \boeta \right) &= Z(H_{e};a,b,c) \\
ac
\left(  \doeta +  \dteta + c \boeta \right)
&=
Z(A_{e}/e;a,b,c).
\end{split}
\end{equation}

We need to be able to determine the polynomials  $\boeta$ ,
$\doeta$ and $\dteta$ uniquely. When $a=1$ we can solve
(\ref{eq:unsolv}) for $c$ and obtain a formula for $Z(\hG;1,b,c)$,
and when $c=1$, we can solve for $a$  to obtain a formula for the
Tutte polynomial $Z(\hG;a,b)$. However, in general the system
(\ref{eq:unsolv}) does not have a unique solution.
 In Section~\ref{sec:plan} we got around this difficulty by
restricting the topology of the ribbon graphs $H_e$. In
Section~\ref{sec:brgeo} we will determine the polynomials  $\boeta$ ,
$\doeta$ and $\dteta$ by considering geometric ribbon graphs. But
before we do this we observe that we could determine these polynomials
uniquely if we used multivariate polynomials $Z(H_e ; a, \vx ,c)$.

\medskip

Label all of the edges of $A_e$ with elements of a set $\vx$ and
consider
the multivariate \BR polynomial $Z(A_e; a,\vx ,c)$. Each state $s$ of
$H_e$ is uniquely determined by the set of edges it contains and
therefore
gives rise to a unique monomial $\prod_{e \in E(s)} x_e$ in $\vx$.
This
in turn determines a unique term of $Z(A_e; a,\vx ,c)$. This means
that
each monomial term in $Z(A_e; a,\vx ,c)$ appears exactly once on
each
side of equation (\ref{eq:fullunsolv}) and we can therefore solve 
(\ref{eq:unsolv}) by comparing terms. Lemma~\ref{lem:zgz} then
gives:
\begin{theorem}
Let $\hG$ be a ribbon graph with the 2-decomposition
$\left(G,\{H_e\}_{e\in E}\right)$, and let $A_e$ be the ribbon graph
$H_e$
with an additional ribbon $e$ joining the vertices $u_e$ and $w_e$.
Then
\[ Z (\hG;a, \vx ,c) = 
\mathcal{G} \left(
Z\left(\tG; a, \vx ,
c  \right)
\right),
\]
where $p_e$, $q_e$ and $r_e$ are uniquely determined by the pair of
equations
\begin{align*}
cp_e+q_e+r_e &=  Z(A_{e}/e;a,\vx,c) \\
p_e+acq_e+cr_e  &=  Z(H_{e};a,\vx,c),
\end{align*}
and $\mathcal{G}$ is induced by
\[
\mathcal{G} : \prod_{e\in E} f_e^{\alpha_e}g_e^{\beta_e} \mapsto
\prod_{e\in E} \left(  \frac{r_{e}}{ac}  \right)^{\alpha_e \beta_e}
\left( \frac{q_e}{2ac} \right)^{(1-\beta_e)}
c^{(\alpha_e \beta_e-\alpha_e)}
\left(\frac{p_e}{ac} \right)^{(\beta_e-\alpha_e \beta_e)} 
.
\]
\hfill$\Box$
\end{theorem}

\subsection{Geometric ribbon graphs} \label{sec:brgeo}

Let $H_e$ be a ribbon graph (so that $t(H_e)=0$). Previously we
considered the ribbon graph $A_e$, which consisted of $H_e$ with an
additional untwisted ribbon $e=(u_e,w_e)$. Here we consider the
geometric ribbon graph $A_{\te}$,
which is obtained from $H_e$ by inserting a half-twisted ribbon $\te$
between the vertices $u_e$ and $w_e$.
This ribbon is inserted into the ribbon graph $H_e$ according to the
conventions of Notation~\ref{notation1}. $A_{\te}$ will always denote
a
geometric ribbon graph constructed from $H_e$ in this way. We will
only need to distinguish the edge $\te$ in $Z(a,\vx,c,d)$ so as
before,
we will denote by $\vb$ the specialization of $\vx$ which sets $x_e =
b$
for $e\neq q$.

The insertion of the ribbon $\te$ into a state $s \in \s (H_e)$
determines a unique state $\ts \in \s (A_{\te})$. Notice that since
$s$ is orientable, $\ts$ is
non-orientable (and $t(\ts)=1$) if and only if the vertices $u_e$ and
$w_e$ lie in the same connected component in $s$.  We will use this
observation in the proof of the following theorem.

\begin{theorem}\label{th:full}
Let $\hG$ be a ribbon graph with the 2-decomposition
$\left(G,\{H_e\}_{e\in E}\right)$, and let $A_{\te}$ be the ribbon
graph $H_e$ with a half-twisted edge $\te=(u_e, w_e)$ inserted.
Then
\[ Z (\hG;a, b ,c) =
\mathcal{G}
\left( 
Z\left(\tG; a, \vx ,
c  \right)\right)
\]
where
\[
q_e = Z\left(A_{\te}/\te;a,b,c,0\right);
\]
$p_e$ and $r_e$ are uniquely determined by the pair of equations 
\begin{align*}
cp_e+r_e &= Z(A_{\te}/\te;a,b,c,1)  -
Z(A_{\te}/\te;a,b,c,0)
\\
p_e+cr_e  &=  Z(H_{e};a,b,c) - ac
Z(A_{\te}/\te;a,b,c,0);
\end{align*}
and $\mathcal{G}$ is induced by
\[
\mathcal{G} : \prod_{e\in E} f_e^{\alpha_e}g_e^{\beta_e} \mapsto
\prod_{e\in E} \left(  \frac{r_e}{ac}  \right)^{\alpha_e \beta_e}
\left( \frac{q_e}{2ac} \right)^{(1-\beta_e)}
c^{(\alpha_e \beta_e-\alpha_e)}
\left(\frac{p_e}{ac} \right)^{(\beta_e-\alpha_e \beta_e)} 
.
\]
\end{theorem}
\begin{proof}
Just as in the derivation of (\ref{eq:unsolv}), we consider the
effect of the insertion of the edge $\te$ into a state of $H_e$ on the
number of connected components, edges, boundary cycles, and the
orientability of the surface. Let $\dcH$,  $\ddH$ and  $\bcH$ be the
partition of $\s (H_e)$ described in Subsection~\ref{ss:brCyDe}.
Suppose that $s \in \s (H_e)$. 
Then if $s \in \dcH$ the insertion of $\te$ reduces the number of
boundary cycles by one and makes the surface non-orientable; if $s \in
\ddH$ the insertion of $\te$ reduces the number of boundary cycles by
one and decreases the number of connected components by one; if $s \in
\bcH$ the insertion of $\te$ makes the surface non-orientable but the
number of boundary cycles and connected components is unchanged.

Now since every state in $\s(A_{\te})$ is either a state in  $\s
(H_e)$
or a state determined by the insertion of the ribbon $\te$ into a
state in
$\s (H_e)$, we have
\begin{multline*}
Z(A_{\te}; a, \vb ,c,t) 
=
(x_{\te} c^{-1}t + 1) \sum_{s \in \dcH}
a^{k(s)}b^{e(s)}c^{\partial(s)} + 
(x_{\te} a^{-1}c^{-1} + 1) \sum_{s \in \ddH}
a^{k(s)}b^{e(s)}c^{\partial(s)}  \\  + 
(x_{\te} t + 1) \sum_{s \in \dcH}
a^{k(s)}b^{e(s)}c^{\partial(s)}.
\end{multline*} 
Rewriting this in terms of the state sums  in  (\ref{eq:ribOGS}) gives
\[ 
Z(A_{\te}; a, \vb ,c,t) 
=
(x_{\te} c^{-1}t + 1) ac^{2}\; \doeta
 +
(x_{\te} a^{-1}c^{-1} + 1) a^{2}c^{2} \; \dteta +
(x_{\te} t + 1)ac \; \boeta.
\]
Then, writing this as a linear equation in $x_{\te}$, we have
\begin{equation}\label{eq:twist}
ac   \left( c \doeta +ac\dteta + \boeta \right)
+
ac
\left(  t\doeta +  \dteta + t \boeta \right)x_{\te}
=
Z(H_{e};a,b,c)+x_{\te}Z(A_{\te}/{\te};a,b,c,t),
\end{equation}
where the right hand side is obtained from the deletion-contraction
relation.
Noting that $x_{\te}ac\dteta = x_{\te}Z(A_{\te}/{\te};a,b,c,0)$,
the theorem then follows easily from Lemma~\ref{lem:zgz}.
\end{proof}

\medskip

As an application of our methods we prove that, just as with the
Tutte polynomial, the \BR polynomial is well defined with respect to
the 2-sum of ribbon graphs. Given two ribbon graphs, each with a
distinguished ribbon, there are four ways of forming the 2-sum,
coming 
from the four ways in which the distinguished ribbons can be
identified. The following result tells us that the \BR polynomial
cannot tell the difference, even if the two resulting ribbon graphs
are non-isomorphic.
\begin{proposition}\label{prop:2sum}
Let $F$ and $F^{\prime}$ be two ribbon graphs which are the 2-sums of
the same pair of ribbon graphs along the same distinguished edges.
Then $R(F) = R(F^{\prime})$.
\end{proposition}
\begin{proof}
Both ribbon graphs have the same 2-decomposition $\left( C_2 , H_1,
H_2\right)$. The difference in the ribbon graph arises from the
choices of which pairs of vertices are identified in the formation of
$F$. The result then follows from Theorem~\ref{th:full} since all of
the
formulae are independent of this choice.
\end{proof}

\begin{remark}
Ideally we would like to be able to extend our results to the
calculation of the \BR polynomial of any geometric ribbon graph,
but the decompositions used here do not lend themselves well to
non-orientable geometric ribbon graphs. The problem is
that if   $\hG$ non-orientable it does not follow that one of $G$ or
the  $H_e$  are non-orientable, therefore the decompositions of
Subsections~\ref{ss:brCyDe} do not record the
orientability of the original states of $\hG$. 
\end{remark}

\begin{remark}
The underlying ideas in this paper extend to $k$-sums (or more
generally $k$-decompositions) of ribbon graphs and we believe that
our results will also extend to $k$-sums of ribbon graphs.  The main
difficulty in generalizing the results appears to be in showing that
the analogues of (\ref{eq:twist}), which arise by considering $H_e
\oplus_k K_k$, have a unique solution.
 
\end{remark}

\section{An application to knot theory}
As mentioned previously, our main motivation for this work came from
recent results connecting the \BR polynomial and knot polynomials
(\cite{CP, CV, Detal, Mo})  which generalize well known relations
between the Tutte polynomial and knot polynomials (\cite{Ja, Th}). In
particular we were interested in generalizing connections between the
behaviour of the Jones polynomial of an alternating link and the
matroid properties of the Tutte polynomial discussed by the first
author in~\cite{Hu}. As an application we will show how invariance of
the Jones polynomial under the mutation of a link can be explained in
terms of the behaviour of ribbon graph polynomials under the 2-sum.
We will assume a familiarity with basic knot theory.

\begin{figure}
\[\epsfig{file=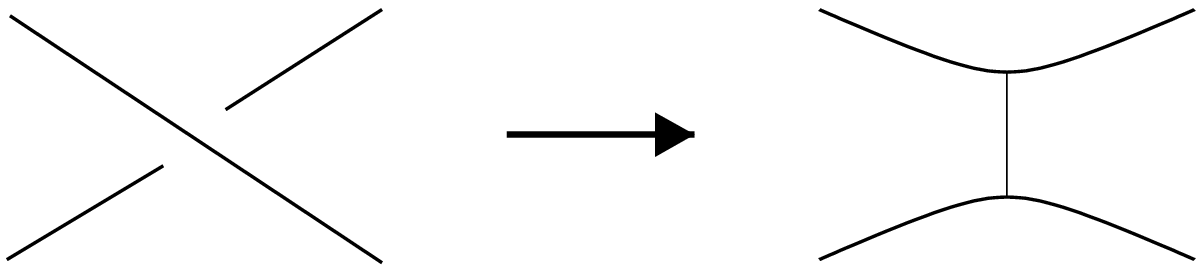, height=1cm}\]
\caption{}
\label{fig:splice}
\end{figure}

\begin{figure}
\begin{center}
\subfigure[]{\epsfig{file=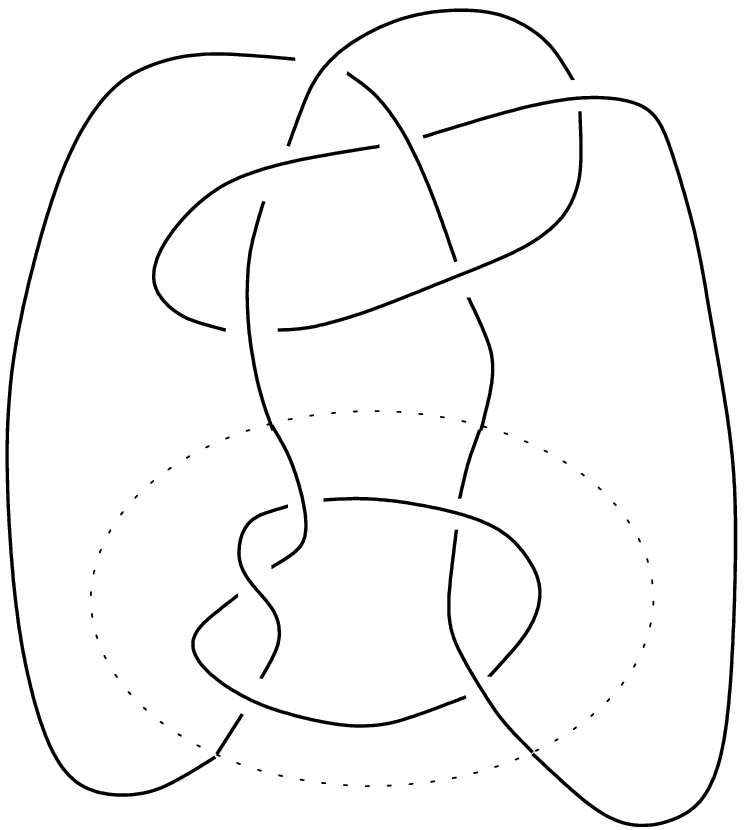, height=3cm}} \quad \quad \quad
\subfigure[]{\epsfig{file=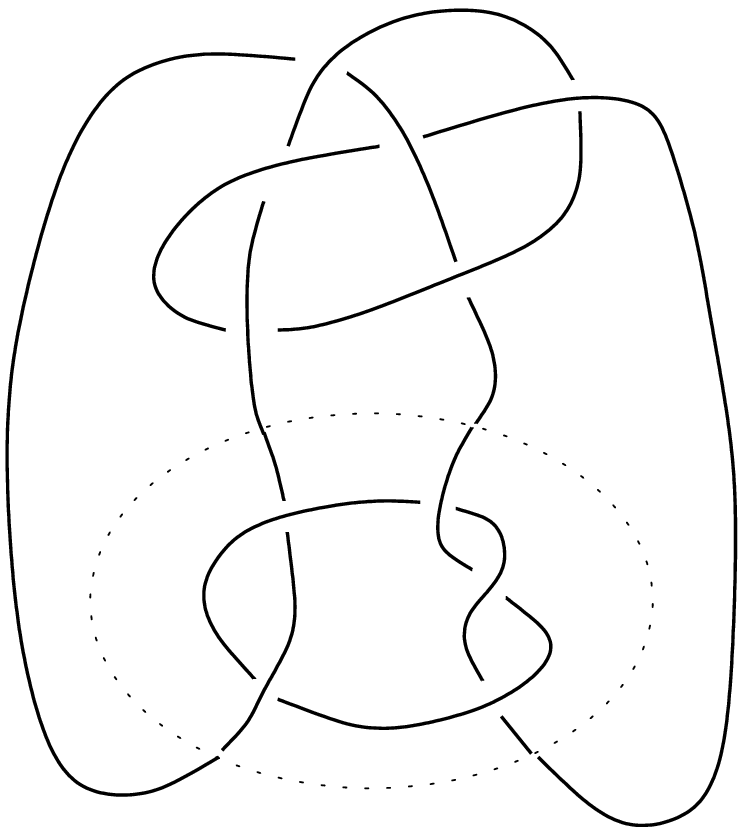, height=3cm}}
\end{center}
\caption{}
\label{fig:mutants}
\end{figure}

\begin{figure}
\begin{center}
\subfigure[]{\epsfig{file=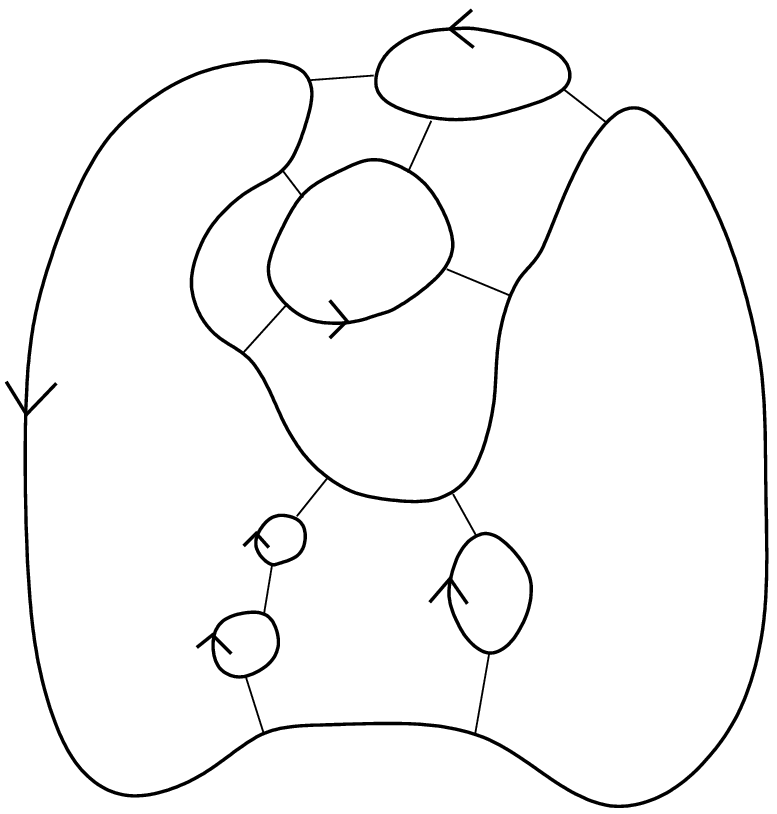, height=3cm}} \quad \quad \quad
\subfigure[]{\epsfig{file=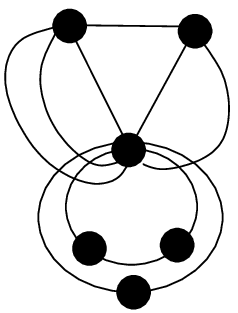, height=3cm}}
\end{center}
\caption{}
\label{fig:ktexamp}
\end{figure}

We will begin by reviewing the construction in \cite{Detal} of a
ribbon graph from a link diagram.
Given a link diagram $D$, begin by replacing each crossing with its
A-splicing as shown in Figure~\ref{fig:splice}.
This gives a collection of disjoint circles in the plane, which we
will call {\em cycles}, and  arcs which record the splicing. As noted
in \cite{Detal},
there is a unique orientation of the cycles in such
a way that the outermost cycles inherit their orientation from the
plane and such that whenever two cycles are nested, they have the
opposite orientation. We will denote the resulting diagram
$\mathcal{D}$. We can then define a ribbon graph $F_D$ by associating
a vertex with each cycle of $\mathcal{D}$ and putting an edge between
two vertices whenever there is an arc in $\mathcal{D}$ connecting the
corresponding cycles. The cyclic ordering of the incident half-edges
at a vertex of $F_D$ is taken to be the cyclic order of the
corresponding arcs on the corresponding cycle in $\mathcal{D}$.
Figure~\ref{fig:ktexamp} shows the diagram $\mathcal{D}$ and the
associated ribbon graph for the Kinoshita-Terasaka knot shown in
Figure~\ref{fig:mutants}(a).

Recall that the Kauffman bracket $\langle D \rangle \in
\mathbb{Z} [A,A^{-1}]$ is a regular isotopy invariant of
links~\cite{Ka}. It is related to the Jones polynomial $J(L)\in
\mathbb{Z} [t,t^{-1}]$ through the identity
\begin{equation}\label{eq:jones} J(L) = \left.\left(-A^{-3\omega (D)}
\langle D \rangle \right)\right|_{A=t^{-1/4} },\end{equation}
where $D$ is any diagram for $L$ and $\omega (D)$ is the {\em writhe}
of $D$. 
 
The following result from \cite{Detal} generalizes to all links
Thistlethwaite's well known result \cite{Th} relating the Tutte
polynomial and the Jones polynomial of an alternating link.
\begin{theorem}\label{thm:brjones}
Let $\langle D \rangle \in \mathbb{Z} [A,A^{-1}]$ be the Kauffman
bracket of a link diagram $D$ and $F$ be its associated ribbon graph.
Then
\[ \langle D \rangle = A^{n(F)-r(F)} R\left(F; \, -A^4 ,\, -1-A^{-4},
\, (-A^2-A^{-2} )^{-1} \right)  . \]
Consequently the Jones polynomial of a link can be obtained as the
evaluation of the \BR polynomial of an associated ribbon graph.
\end{theorem}

\medskip

Two link diagrams $D$ and $D^{\prime}$ are said to be {\em mutants}
if there exists a circle $C$ in the plane (regarded as $z=0$ in
$\mathbb{R}^3$) which intersects $D$ transversally  in exactly four
points such that by rotating the interior $C$ by $\pi$ radians about
the $x$-, $y$- or $z$-axis gives the diagram $D^{\prime}$. We say two
links are mutants if they admit mutant diagrams. 
Figure~\ref{fig:mutants} shows the Kinoshita-Terasaka (a) and Conway
(b)
knots, which are perhaps the most famous examples of mutant knots. 

We use the results of Section~\ref{sec:brII} to provide a new
perspective on the following well known result.
\begin{proposition}
Let $L$ and $L^{\prime}$ be mutant links admitting mutant diagrams
$D$ and $D^{\prime}$. Then $\langle D \rangle = \langle D^{\prime}
\rangle $.
Moreover when the mutation operation on $D$ respects the orientation
of the diagrams, $J(L)=J(L^{\prime})$.
\end{proposition}   
\begin{proof}
We will show that the ribbon graphs associated with $D$ and
$D^{\prime}$ are 2-sums of the same pair of ribbon graphs (which can
differ due to the ambiguity in the 2-sum). The proposition will then
follow by Proposition~\ref{prop:2sum}, Theorem~\ref{thm:brjones} and
Equation~\ref{eq:jones}.

Since the Kauffman bracket is a  regular isotopy invariant, we may
assume that $D$ is of the form shown schematically in
figure~\ref{fig:mutproof}(a).
   \begin{figure}
\begin{center}
\subfigure[]{\epsfig{file=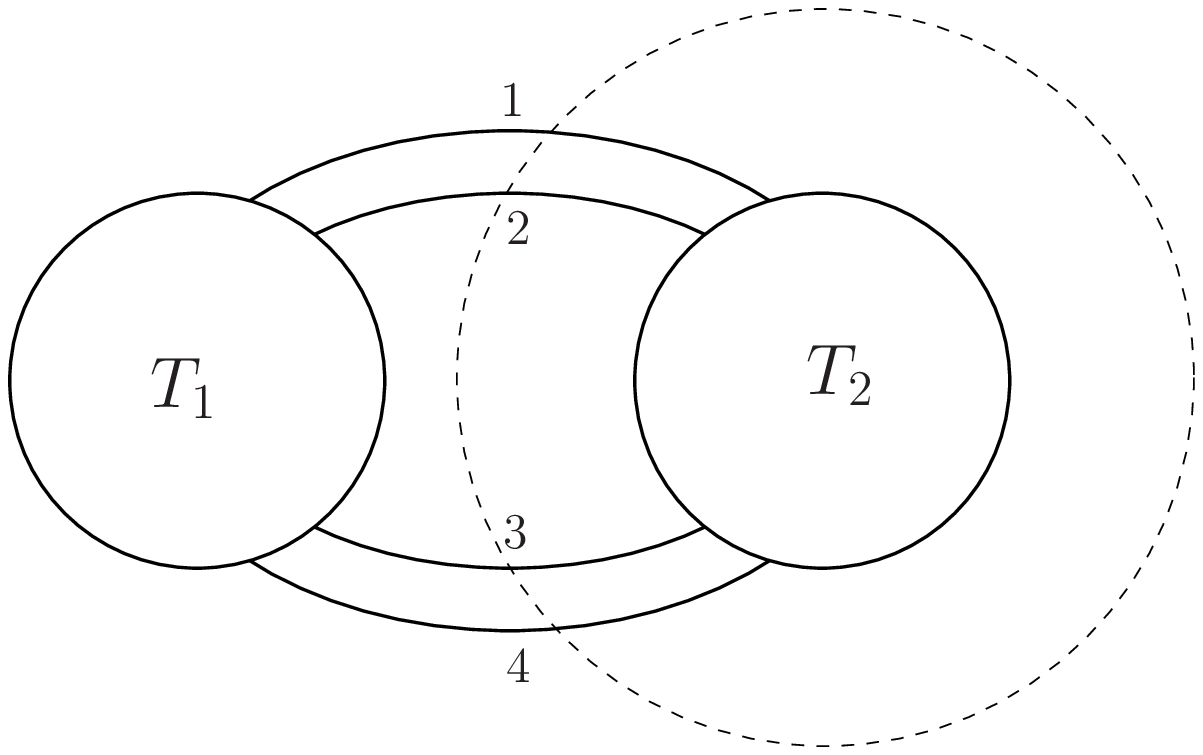, height=3cm}} \quad \quad \quad
\subfigure[]{\epsfig{file=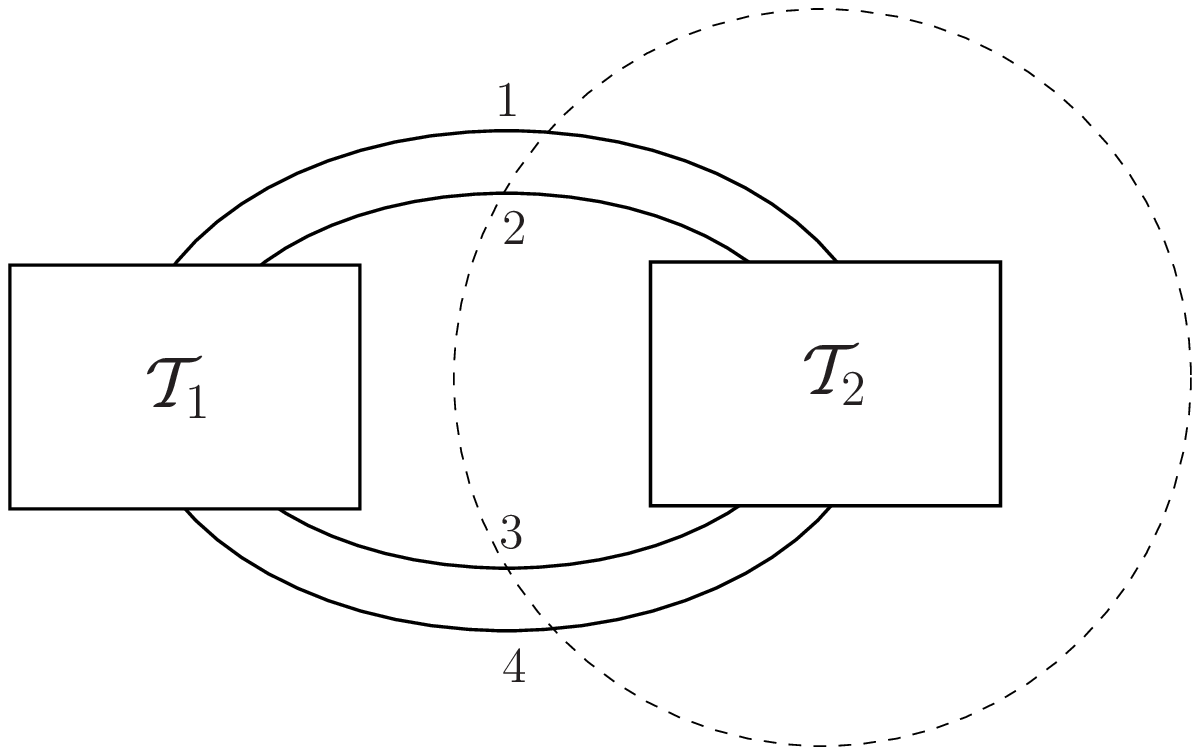, height=3cm}}
\end{center}
\caption{}
\label{fig:mutproof}
\end{figure}
In the figure, $T_{1}$ and $T_{2}$ denote some tangles and the
dotted circle is the circle $C$ used in the definition of mutation.
Consider the intersection points of $C$ and $D$ marked 1, 2, 3, and
4 in the figure. 
These points will lie on cycles in the diagram $\mathcal{D}$. 
Also, since $T_{1}$ and $T_{2}$ share no crossings, cutting
$\mathcal{D}$ at these four points will disconnect the diagram
$\mathcal{D}$. Therefore $\mathcal{D}$ can be represented
schematically as in Figure~\ref{fig:mutproof}(b), where all of the
edges and other cycles are contained in the two boxes.  

Since all of the cycles in $\mathcal{D}$ are closed, there exists a
cycle which contains the pairs of points 1 and 2, 1 and 3, or 1 and
4.

If a cycle contains both the points 1 and 3, then since the cycles of
$\mathcal{D}$ are disjoint, it must contain all of the points 1, 2, 3,
and 4.  This can be done in two ways: with cyclic order (1,2,3,4) or (1,4,3,2).
 Thus there are four cases to consider: when there is a
cycle containing 1 and 2 but not 3 and 4;  1 and 4 but not 2 and 3;
and the two cases when 1, 2, 3, and 4 are in the same cycle. 

If a cycle contains 1 and 2 then there is also a cycle containing 3
and 4. Then, as is clear from
Figure~\ref{fig:mutproof}(a), there are vertices of the ribbon graph
with cyclic ordering of the form
$( 1,  \text{edges of } \mathcal{T}_{1}, 2,  \text{edges of }
\mathcal{T}_{2})$ and  $( 3, \text{edges of }\mathcal{T}_{1}, 4,
\text{edges of }\mathcal{T}_{2})$ (or the inverse orders) and
thus by cutting the discs of the ribbon graph along the  interior
arcs (1,2) and (3,4) we see that the ribbon graph is a 2-sum.     

If a cycle contains  1 and 4 and another cycle contains 2 and 3 then
in the ribbon graph there is a disc with cyclic ordering of the form
$( 1, \text{edges of }\mathcal{T}_{1}, 4,  \text{edges of
}\mathcal{T}_{2})$ and  $( 2,  \text{edges of } \mathcal{T}_{1}, 3,
\text{edges of } \mathcal{T}_{2})$ (or the inverse orders)  and
thus by cutting the discs of the ribbon graph along the  interior
arcs (1,4) and (2,3) we see that the ribbon graph is a 2-sum.     

Similarly, when 1, 2, 3, and 4 are all points on the same cycle (with either order) then
there is a single disc in the ribbon graph such that cutting along
interior arcs (1,2) and (3,4) will separate the ribbon graph into its
2-sum components. 

Finally, consider the rotation of the interior of $C$ on $D$. The
corresponding effect on the diagram $\mathcal{D}$ is the same
rotation of the interior of $C$. It is clear that in the ribbon
graph, this corresponds to changing the way we form the 2-sum of the
two ribbon graphs obtained above.      

\end{proof}

\end{document}